\documentclass{amsart}
\usepackage{amsrefs}
\usepackage{bbm}
\usepackage{accents}
\usepackage{amsmath}
\usepackage{algorithm}
\usepackage{algorithmic}
\usepackage{amsthm}
\usepackage{amssymb} 
\usepackage{graphicx}
\usepackage{epstopdf}
\usepackage{epsfig}
\usepackage{ifthen}
\usepackage{lipsum}
\usepackage{mathrsfs}	
\usepackage{multimedia}
\usepackage{multirow}
\usepackage{wrapfig}
\usepackage{mathtools}
\usepackage{cleveref}
\usepackage{xcolor}
\usepackage{makecell}
\usepackage{cases}
\usepackage[overload]{empheq}

\usepackage[caption=false]{subfig}

\newtheorem{theorem}{Theorem}[section]

\newtheorem{Lem}[theorem]{Lemma}

\theoremstyle{definition}
\theoremstyle{remark}

\newtheorem{Def}[theorem]{Definition}
\newtheorem{Exam}[theorem]{Example}

\numberwithin{equation}{section}

\newcommand{\calu}{\mathcal{U}}
\DeclarePairedDelimiter\ceil{\lceil}{\rceil}
\DeclarePairedDelimiter\floor{\lfloor}{\rfloor}

\DeclareMathOperator{\interior}{int}
\newcommand{\norm}[1]{\left\lVert#1\right\rVert}

\newlength{\dhatheight}

\newlength{\figwidth}
\begin{document} 

\title{Hausdorff Moment Transforms and Their Performance}

\author{Xinyun Wang}
\address{Dept. of Electrical Engineering, University of Notre Dame, Notre Dame, IN 46556}
\email{xwang54@nd.edu}
\thanks{This work was supported in part by the US National Science Foundation through Grant 2007498. }

\author{Martin Haenggi}
\address{Depts. of Electrical Engineering and Applied and Computational Math. and Statistics, University of Notre Dame, Notre Dame, IN 46556}
\email{mhaenggi@nd.edu} 

\begin{abstract}
	Various methods have been proposed to approximate a solution to the truncated Hausdorff moment problem. In this paper, we establish a method of comparison for the performance of the approximations. Three ways of producing random moment sequences are discussed and applied. Also, some of the approximations have been rewritten as linear transforms, and detailed accuracy requirements are analyzed. Our finding shows that the performances of the approximations differ significantly in their convergence, accuracy, and numerical complexity and that the decay order of the moment sequence strongly affects the accuracy requirement.
\end{abstract}

\maketitle

%

\section{Introduction} 

\subsection{Problem formulation and related work}

In many applications, such as computed tomography (CT) image reconstructions \cite{wang2001image} and meta distribution reconstructions \cite{haenggi2016meta}, distributions of bounded support need to be reconstructed from given moment sequences. Without loss of generality, we assume the support of the distributions is $[0,1]$. The problem can be formulated as follows. Let $[n] \triangleq \{1,2,...,n\}$ and $[n]_0 \triangleq \{0\} \cup [n]$. 
Given a finite sequence $\mathbf{m} \triangleq (m_k)_{k=0}^{n}, \, n \in \mathbb{N}$, find/approximate an $F$ such that
\begin{align}
	\int_{0}^{1} x^k \;dF(x) = m_k, \quad \forall k \in [n]_0, \label{eq:truncated}
\end{align} 
where $F$ is right-continuous and increasing with $F(0^-) = 0$ and $F(1) = 1$, i.e., $F$ is a cumulative distribution function (cdf). This problem is known as the truncated Hausdorff moment problem (THMP) \cite{shohat1950problem}. We call a method that finds or approximates a solution $F$ a \emph{Hausdorff moment transform}. In the following, we assume the existence of solutions to the THMP, and $F$ is the complementary distribution function (cdf) of the corresponding distribution. 

A variety of such transforms have been proposed in the literature: \cite{tagliani1999hausdorff} applies the maximum entropy (ME) method \cite{kapur1989maximum} to find a distribution that exactly corresponds to the finite sequence. \cite{mnatsakanov2003some} uses a mixed binomial cdf to approximate such a distribution, which we call the binomial mixture (BM) method. \cite{shohat1950problem}, \cite{guruacharya2018approximation}, and \cite{schoenberg1973remark} propose different approximations based on the Fourier-Jacobi (FJ) series.  Besides, Markov derives the (pointwise) suprema and infima of the cdfs among all solutions of the THMP \cite{shohat1950problem}. If the integer moments can be extended to the complex plane, the Gil-Pelaez (GP) theorem \cite{gil1951note} can be applied, henceforth called the GP method.\footnote{Other forms of  Fourier inversion  may also be applied, such as the inverse Mellin  transform \cite{penson2022exact}.} 

Each of the existing methods has its limitations. The GP method yields an integral of infinite range that is too complicated to be analytically calculated. Its numerical approximation requires a careful selection of integration upper limits and step sizes, and the integrand may contain a singularity. However, even if accurate results are obtained, the GP method does not solve \eqref{eq:truncated} since it is based on characteristic functions (or imaginary moments), and it is not always possible to extend integer moments to the complex plane. The ME method is not useful for large $n$ due to the ill-conditioning of the Hankel matrices used in the calculation \cite{NOVIINVERARDI2003frac}. The BM method usually requires a large number of integer moments to reach sufficient accuracy. {The FJ approximations used in \cite{shohat1950problem}, \cite{guruacharya2018approximation}, and \cite{schoenberg1973remark} are functional approximations which may produce unreasonable results that violate the general properties of probability distributions. Further, a detailed comparison and discussion of convergence and concrete expressions for the approximations is missing in \cite{shohat1950problem}, \cite{guruacharya2018approximation}, and \cite{schoenberg1973remark}.} 
Hence, there is a need to study, improve, and compare the existing methods in detail and develop new ones with provably good numerical performance.

\subsection{Contributions}

\begin{itemize} 
	\item We establish a method of comparison for the performance of the HMTs, based on their accuracy, computational complexity, and accuracy requirements. This way, the most suitable HMT for given objectives can be chosen. The key point to make fair comparisons is to generate representative moment sequences. The problem of 
	generating random moment sequences is non-trivial since the probability that sequences i.i.d.~ on $[0,1]^n$ are moment sequences is about $2^{-n^2}$ \cite[Eq.~1.2]{chang1993normal}. We provide three different ways to solve it: One is to generate random integer moment sequences that are uniformly distributed in the moment space; the other two are based on analytic expressions, either of moments or of the cdfs.
	\item To find or approximate a solution to the THMP, we propose a method based on the average of the infima and suprema from the Chebyshev-Markov (CM) inequalities and a method based on the Fourier-Chebyshev (FC) series, called the FC method. We also provide the concrete expressions of the FJ approximation put forth in \cite{shohat1950problem}. As some of the HMTs are essentially functional approximations, we propose a mapping to restore properties of probability distributions. 
	\item We rewrite some of the HMTs as linear transforms of the given sequence. The transform matrix depends only on the length of the sequence. It can be calculated offline to (further) save time. Besides, the linear transformation allows for a direct analysis of the accuracy requirements. 
\end{itemize}

\section{Preliminaries}
\subsection{Existence and uniqueness of solutions}
Let $\mathcal{F}_n$ be the family of $F$ that solves \eqref{eq:truncated}. The existence and uniqueness of solutions can be verified by the Hankel determinants, defined by
\begin{align}
	\underline{H}_{l} & \triangleq
	\begin{cases}
		\begin{vmatrix}
			m_0 &  \dots & m_{l/2} \\  
			\vdots  & \ddots & \vdots  \\ 
			m_{l/2} & \dots & m_{l}   
		\end{vmatrix}, \quad\text{even } l,\\
	\vspace{0.005\textheight}\\
		\begin{vmatrix}
			m_1 &  \dots & m_{(l+1)/2} \\  
			\vdots  & \ddots & \vdots  \\ 
			m_{(l+1)/2} & \dots & m_{l}   
		\end{vmatrix}, \quad\text{odd }l,
	\end{cases} 
\end{align}
\begin{align}
	\overline{H}_l & \triangleq
	\begin{cases}
		\begin{vmatrix}
			m_1 - m_2 & \dots & m_{l/2} - m_{l/2+1}\\  
			\vdots  & \ddots & \vdots  \\ 
			m_{l/2} - m_{l/2+1} &  \dots &  m_{l-1}-m_{l}   
		\end{vmatrix}, \quad\text{even } l,\\
	\vspace{0.005\textheight}\\
		\begin{vmatrix}
			m_0- m_1 & \dots & m_{(l-1)/2} - m_{(l+1)/2}\\  
			\vdots  & \ddots & \vdots  \\ 
			m_{(l-1)/2} - m_{(l+1)/2} &  \dots &  m_{l-1}-m_{l}   
		\end{vmatrix}, \quad\text{ odd }l.
	\end{cases} 
\end{align}
A solution exists if $\underline{H}_{l}$ and $\overline{H}_l$ are non-negative for all $l \in [n]$. Furthermore, the solution is unique, i.e., $\vert \mathcal{ F}_n \vert = 1$, if and only if there exists an $l \in [n]$ such that $ \underline{ {H}}_l = 0$ or $ \overline{H}_{l} = 0$. If the solution is unique, $ \underline{ {H}}_k =  \overline{H}_{k} = 0$ for all $k >l$, and the distribution is discrete and uniquely determined by the first $l$ moments. \Cref{ex:hankel} shows an example with $n = 2$. More interesting and relevant is the case where the Hankel determinants $ \underline{ {H}}_k$ and $ \overline{H}_{k}$ are positive for all $l \in [n]$ and the problem has infinitely many solutions, i.e.,  $\vert \mathcal{ F}_n \vert = \infty$. In the following, we focus on this case by assuming the positiveness of all the Hankel determinants.

 \begin{Exam}[Two sequences with $n = 2$]\label{ex:hankel}
 	{\leavevmode
 		\begin{enumerate}
 			\item For $\mathbf{m}= (1, 1/2, 1/2)$, we have
 			\begin{align*}
 				\underline{ {H}}_{1} = m_1 = 1/2 > 0,\quad\overline{H}_{1} = m_0 - m_1 = 1/2>0,\\
 				\underline{ {H}}_{2} = \begin{vmatrix}
 					m_0 & m_1\\
 					m_1 & m_2
 				\end{vmatrix} = 1/4 >0,\quad\overline{H}_{2} = 
 				m_1 - m_2 = 0.
 			\end{align*}
 			The uniquely determined discrete distribution has only two jumps, one at $0$ and one at $1$. Hence the probability density function (pdf) is $\frac{1}{2} \delta(x) + \frac{1}{2} \delta(x-1)$.
 			\item For $\mathbf{m} = (1, 1/2, 1/3)$, we have
 			\begin{align*}
 				\underline{ {H}}_{1}  = m_1= 1/2 > 0,\quad\overline{H}_{1}= m_0 - m_1= 1/2  >0,\\
 				\underline{ {H}}_{2} = \begin{vmatrix}
 					m_0 & m_1\\
 					m_1 & m_2
 				\end{vmatrix} = 1/12 >0,\quad\overline{H}_{2} = 
 				m_1 - m_2 = 1/6 >0.
 			\end{align*}
 			There are infinitely many solutions since all Hankel determinants are positive, such as $f(x) = 1 $ and $f(x) = \frac{1}{2} \delta\left(x- \left(3+\sqrt{3}\right)/6 \right) + \frac{1}{2} \delta\left(x- \left(3-\sqrt{3}\right)/6 \right) $. 
 		\end{enumerate}
 	}
 \end{Exam}

\subsection{Adjustment}

Here, we introduce a mapping to restore properties of probability distributions.

Some approximations are essentially functional approximations which may not preserve the properties of cdfs, such as monotonicity and ranges restricted to $[0,1]$. An example can be found in \Cref{fig:FJcon20}. Besides, for practical purpose, the numerical evaluation of approximations are necessarily to be restricted to a finite set of $[0,1]$. Here, we consider a finite set $\mathcal{U}_n$ with cardinality $(n+1)$, and let $(u_0, u_1, u_2, ..., u_{n})$ be the sequence of all elements in $\mathcal{ U}_n$ in ascending order, where $u_0 = 0$ and $u_n = 1$. If there is no further specification, we assume $\mathcal{U}_{n}$ is the uniform set $\{\frac{i}{n}, i \in [n]_0\}$. \Cref{def:correction} introduces a tweaking mapping to adjust the values of the approximations so that they adhere to the properties of cdfs after adjustments. 
 
 \begin{Def}[Tweaking mapping]\label{def:correction} 
 	We define $T$ that maps $F: \mathcal{ U}_n\mapsto \mathbb{R}$ to $\hat{F}: \mathcal{ U}_n \mapsto [0,1]$ such that  
 	\begin{enumerate}
 		\item $T\left(F\right)(u_i) =\min\left(1, \max\left(0, F(u_i)\right)\right)$, $\forall i \in [n]_0$.
 		\item $T\left(F\right)(0) = 0 $ and $T\left(F\right)(1) = 1$. 
 		\item $T\left(F\right)(u_{i+1})  = \max(F(u_i), F(u_{i+1}))$, $\forall i \in [n-1]_0$. 
 	\end{enumerate} 
 \end{Def}
 
 The tweaking mapping in \Cref{def:correction} eliminates the out-of-range points and keeps the partial monotonicity. Let $F:\calu_{n} \mapsto \mathbb{R}$ be an approximation of a solution to the truncated HMP. After tweaking, we obtain monotone and bounded approximations over a discrete set, denoted as $\hat{F} = T\left(F\right)$. Next, we interpolate the discrete points to obtain continuous functions. Specifically, we apply monotone cubic interpolation to $\hat{F}$. Monotone cubic interpolation preserves monotonicity.  For convenience, we refer to the function after interpolation as $\hat{F}$ as well, and we call the original function $F$ the \emph{raw} function and the tweaked and interpolated function $\hat{F}$ the \emph{polished} function. We use $\hat{F}|_{\calu_n}$ to emphasize that the polished function $\hat{F}$ is originally sampled over $\calu_n$. The polished function is continuous, monotone and $[0,1] \mapsto [0,1]$. 
 
 \subsection{Distance metrics}
In this part, we introduce two distance metrics between two  functions $F_1, F_2: [0,1] \mapsto [0,1]$ are.

\begin{Def}[Maximum and total distance] 
	The maximum distance between functions $F_1$ and $F_2$ is defined as the  $\mathcal{L}_{\infty}\mbox{-}\text{norm}$ of $F_1$ and $F_2$ over $[0,1]$, i.e., 
	\begin{align*}
		d_{\rm{M}}(F_1, F_2)  & \triangleq  \norm{F_1 - F_2}_{\mathcal{L}_{\infty}[0,1]},\\
		& = 	\max_{x \in [0,1]} \lvert F_1(x) - F_2(x)  \rvert.
	\end{align*}
The total distance between functions $F_1$ and $F_2$ is defined as the $\mathcal{L}_{1}\mbox{-}\text{norm}$ of $F_1$ and $F_2$ over $[0,1]$, i.e., 
	\begin{align*}
		d(F_1, F_2)  & \triangleq  \norm{F_1 - F_2}_{\mathcal{L}_{1}[0,1]},\\
		& = \int_0^1 \lvert F_1(x) - F_2(x)  \rvert dx.
	\end{align*} 
\end{Def} 

We mainly focus on the total distance and use the maximum distance as a secondary metric. As the maximum and total distance are designed to be applied to the polished functions, the exact value of the total distance and maximum distance of the linear/monotone cubic interpolation of the polished functions can be calculated.

\section{Solutions and approximations}
\subsection{The CM inequalities}\label{sec:cm}
Markov \cite{shohat1950problem} provides a method to obtain $\inf_{F\in \mathcal{F}_n} F(x_0)$ and $\sup_{F\in \mathcal{F}_n} F(x_0)$ for any $x_0 \in [0,1]$. The most important step of the method is the  construction of a discrete distribution in $\mathcal{ F}_n$ where the maximum possible mass is concentrated at $x_0$. The details of the method can be found in \cite{markoff1884demonstration} and \cite{shohat1950problem}, and explicit expressions can be found in \cite{zelen1954bounds} and \cite{wang2022cm}. The inequalities established by the infima and suprema obtained in \cite[Thm.~1]{wang2022cm} are also known as the CM inequalities.
In contrast, the infima and suprema given by \cite{Rcz2006AMB} are for distributions where the  support can be an arbitrary Borel subset of $\mathbb{R}$. As a result, they are loose for distributions with support $[0,1]$.

\begin{figure}
	\centering
	\subfloat[$m_k = 1/(k+1),\quad  k \in {[n]}_0$.]{\includegraphics[width =0.48 \figwidth]{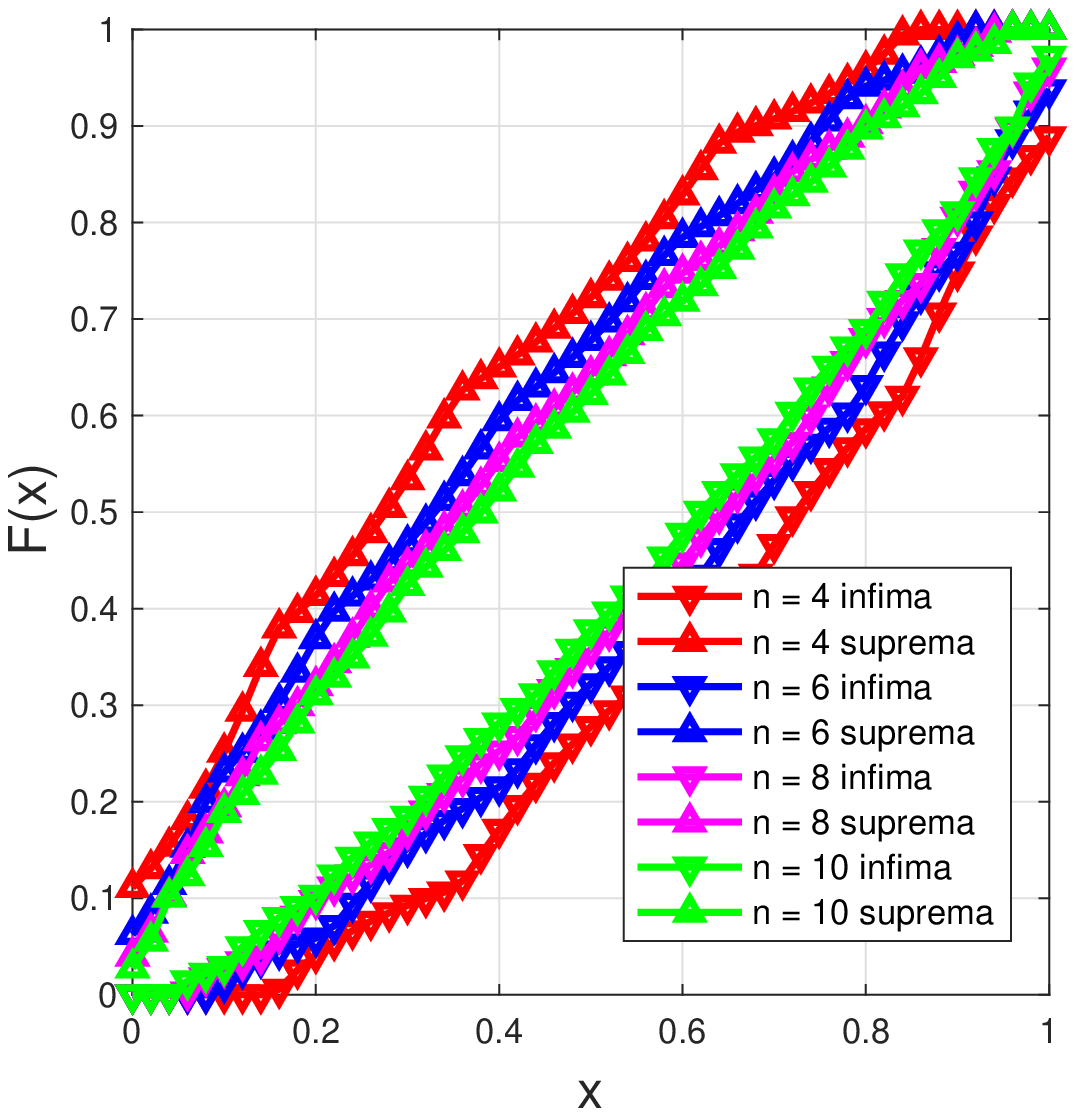}}\quad
	\subfloat[$m_k = \frac{1}{5}( 8^{-k} + 3^{-k} + 2^{-k} + \left(\frac{3}{2}\right)^{-k} + \left(\frac{5}{4}\right)^{-k} 
	),\quad k \in {[n]}_0$.\label{fig:discrete}]{\includegraphics[width =0.48 \figwidth]{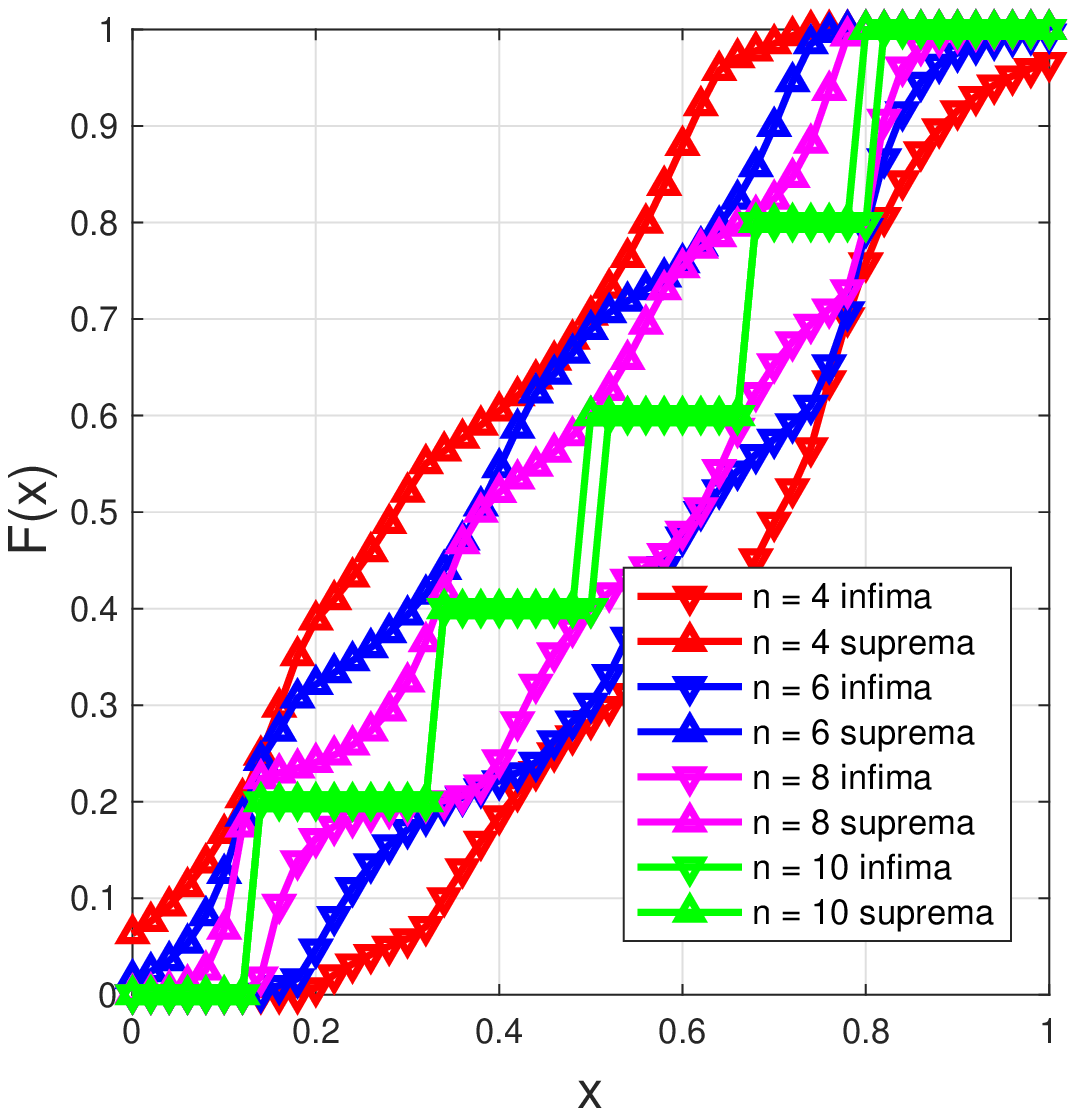}} 
	\caption{ The infima and suprema from the CM inequalities for $n = 4, 6, 8, 10$. $x$ is discretized to $\mathcal{U}_{50} = \{ i/50, \, i \in [50]_0\}$.}  \label{fig:markovbounds}
\end{figure}

\Cref{fig:markovbounds} shows two examples of the infima and suprema for different $n$. The convergence of the infima and suprema as $n$ increases is apparent.  The discontinuities of the derivatives in \Cref{fig:markovbounds} coincide with the roots of the orthogonal polynomials in \cite[Thm.~1]{wang2022cm}, as stated in the second part of \cite[Lem.~1]{wang2022cm}. {For \Cref{fig:discrete}, $\underline{H}_{9} = 0$ indicates a discrete distribution.} The infima and suprema get drastically closer from $n = 8$ to $n = 10$ because $n = 10$ is sufficient to characterize a discrete distribution with $5$ jumps while $n = 8$ is not. {In fact, the infima and suprema match almost everywhere (except at the jumps) for $n = 10$, but \Cref{fig:discrete} only reveals the mismatch at $0.5$ and $0.8$ since only $0.5$ and $0.8$ are jumps that fall in the sampling set $\mathcal{ U}_{50}=\{ i/50, \, i \in [50]_0\}$.  }

\begin{figure}[!htbp]
	\centering
	\includegraphics[width = 0.4\figwidth]{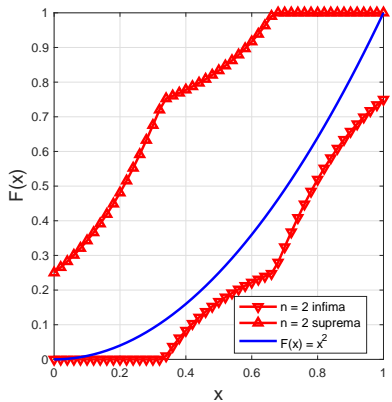}
	\caption{ The infima and suprema from the CM inequalities for $m_n = 1/(n+1)$ and $F(x) = x^2$.}\label{fig:cmsexm}
\end{figure}

Although all solutions to the THMP lie in the band defined by the infima and suprema, the converse is not necessarily true. \Cref{fig:cmsexm} gives an example, where $F$ lies in the band but is not a solution to the THMP. Its corresponding $m_1 $ and $m_2$ are $2/3$ and $1/2$, respectively, instead of $1/2$ and $1/3$ that define the band. 

In the following  subsections, we propose a method based on the FC series and a CM method based on the average of the infima and suprema from the CM inequalities and compare them against the existing methods. We also improve and concretize some of the existing methods. {In particular, we provide a detailed discussion of the existing FJ approximations and fully develop the one mentioned in \cite{shohat1950problem}. The GP method is not applicable to the THMP as it utilizes characteristic functions or imaginary moments instead of integer moments. Besides, the ME method by design yields solutions to the THMP, while the other methods by design obtain approximations of a solution to the THMP.

\subsection{Well-established methods}
\subsubsection{BM method}

A piecewise approximation of the cdf based on binomial mixtures is proposed in \cite{mnatsakanov2003some}. For any positive integer $n$, the approximation\footnote{In this section, the term ``approximations'' is used when the outputs are not always solutions to the THMP.} by the BM method  is defined as follows. 
\begin{Def}[Approximation by the BM method]
	\begin{align}
		{F}_{\mathrm{BM},n}(x)  \triangleq
		\begin{cases}
			\sum_{k = 0}^{\lfloor nx \rfloor} h_k,&\quad x \in (0,1],\\
			0, &\quad x = 0, 
		\end{cases}
	\end{align}
	where $h_k = \sum_{i= k }^{n} \binom{n}{i} \binom{i}{k} (-1)^{i-k} m_i$. 
\end{Def}

Note that $ {F}_{\mathrm{BM},n} \to F $ as $n \to \infty$ for each $x$ at which ${F}$ is continuous \cite{mnatsakanov2003some,mna2008hausdorff}.

\begin{figure}
	\centering
	\includegraphics[width = 0.45\figwidth]{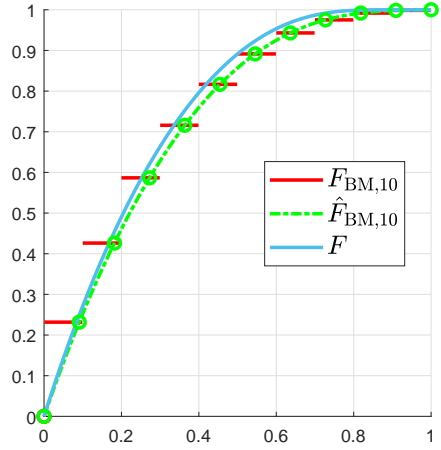}
	\caption{$F_{\mathrm{BM},10}$, $\hat{F}_{\mathrm{BM},10}$ and $F$, where $\bar{F}  = \exp(-\frac{x}{1-x}) (1-x)/ (1-\log(1-x))$. The $\circ$ denote  ${F}_{\mathrm{BM},10}|_{\mathcal{U}_{11}}$.\label{fig:BM_floor}}
\end{figure}		

Let $\hat{F}_{\mathrm{BM},n}$ denote the monotone cubic interpolation of ${{F}_{\mathrm{BM},n}}|_{\mathcal{ U}_{n+1}}$ \cite{haenggi2018efficient}.  As shown in \Cref{fig:BM_floor}, $ {F}_{\mathrm{BM},n} $ is piecewise constant, and for continuous $F$, $\hat{F}_{\mathrm{BM},n}$ provides a better approximation than ${F}_{\mathrm{BM},n}$. Note that ${{F}_{\mathrm{BM},n}}$ is a solution to the THMP only when $n = 1$ and $\hat{F}_{\mathrm{BM},n}$ may not be a solution for any positive integer $n$.

\subsubsection{ME method}

The ME method \cite{tagliani1999hausdorff,kapur1989maximum} is also widely used to solve the THMP. For any positive integer $n$, the solution by the ME method is defined as follows. 

\begin{Def}[Solution by the ME method]
\begin{align}
	f_{\mathrm{ME}, n}(x) \triangleq \exp\left(-\sum_{k = 0}^n \xi_k x^k\right), 
\end{align}
where the coefficients $\xi_k,i \in [n]_0$, are obtained by solving the equations 
\begin{align}
	\int_0^1 x^i \exp\left(-\sum_{k = 0}^n \xi_k x^k\right) \,dx = m_i, \quad i \in [n]_0.\label{eq:me}
\end{align}
\end{Def}

{A unique solution of \eqref{eq:me} exists if and only if the Hankel determinants of the moment sequence are all positive. It can be obtained by minimizing the convex function \cite{tagliani1999hausdorff} 
	\begin{align}
		\Gamma(\xi_1, \dots, \xi_n) \triangleq \sum_{ k = 1}^n \xi_k m_k + \log\left( \int_{0}^{1} \exp\left(-\sum_{k =1}^n \xi_k x^k\right)\right) \label{eq:me_cvx}
	\end{align}
	and solving 
	\begin{align}
		\int_0^1  \exp\left(-\sum_{k = 0}^n \xi_k x^k\right) \,dx = m_0 = 1
	\end{align}
for $\xi_0$.
}
The corresponding cdf is given by 
\begin{align}
	F_{\mathrm{ME}, n}(x) = \int_0^x \exp\left(-\sum_{k = 0}^n \xi_k v^k\right)\, dv.\label{eq:mecdf}
\end{align}

{ Since \eqref{eq:me_cvx} is convex, it can in principle be solved by standard techniques such as the interior point method. As the cdf given by \eqref{eq:mecdf} is a solution to the THMP, its error can be bounded by the infima and suprema from the CM inequalities, which, however, requires additional computation. 
	\cite{tagliani1999hausdorff} shows that $F_{\mathrm{ME}, n} \to F$ in entropy and thus in $\mathcal{L}_1$ as $n \to \infty$.}  But as shown in \cite{NOVIINVERARDI2003frac}, due to ill-conditioning of the Hankel matrices, the ME method is not useful when the number of given moments is large, i.e., we may not be able to solve \eqref{eq:me_cvx} practically.  
  
\subsection{Improvements and concretization of existing methods}\label{sec:imp}

\subsubsection{GP method}

By the Gil-Pelaez theorem \cite{gil1951note}, ${F} $ at which ${F} $ is continuous\footnote{If $F$ is not continuous at some point $x_0$, the result obtained by the Gil-Pelaez theorem is the average of $F(x_0)$ and $F(x_0^-)$.} can be obtained from the imaginary moments by
\begin{align}
	{F} (x) & = \frac{1}{2} - \frac{1}{\pi} \int_0^\infty r(s,x) ds \label{eq:sdf_cdf}, \quad x \in [0,1],
\end{align}
where $r(s,x) \triangleq {\Im{(e^{-js\log x} m_{js})}}/{s}$, $j \triangleq  \sqrt{-1}$, and $\Im(z)$ is the imaginary part of the complex number $z$.  

Since the integral in \eqref{eq:sdf_cdf} cannot be analytically calculated, only numerical approximations can be obtained. An accurate evaluation is exacerbated by the facts that the integration interval is infinite and that the integrand may diverge at $0$. We define the approximate value obtained through the trapezoidal rule as follows. 

\begin{Def}[Approximation by the Gil-Pelaez method]
	\begin{align}
		{F}_{\rm{GP}(\Delta s,\upsilon )}(x) \triangleq 
		\frac{1}{2} - \frac{\Delta s}{2\pi} \sum_{i = 1}^{N} \left ( r(s_i,x) + r(s_{i-1},x) \right ), \quad x \in [0,1],
	\end{align}
	where $N \triangleq  \floor{\frac{\upsilon }{\Delta s}}$ and $s_i \triangleq  i \Delta s$.
\end{Def}
Note that as $\Delta s \to 0$ and $\upsilon \to \infty$, $	{F}_{\rm{GP,\Delta s,\upsilon }} \to {F} $ for each $x$ at which ${F}$ is continuous.

We evaluate ${F}_{\rm{GP}(\Delta s,\upsilon )}$ over the set ${\mathcal{U}_{100}}$ for different $\Delta s$ and $\upsilon$. \Cref{fig:gp_precision} shows the average of the total and maximum  distance of ${F} $ and {polished $\hat{F}_{\rm{GP}(\Delta s,\upsilon )}$} as the step size $\Delta s$ and upper integration limit $\upsilon$ change.  As the step size $\Delta s $ decreases to $0.1$ and as the upper integration limit $\upsilon$ increases to $1000$, the average of the total distance decreases to less than $10^{-3}$, which indicates that  the approximate cdf obtained by numerically evaluating the integral of Gil-Pelaez method is very close to the exact cdf {and thus can be used as a reference}. {Both the integration upper limit and step size are critical in evaluating the integral.} {A general rule for selecting the upper integration limit is to examine the magnitude of $m_{js}$. For example, if $|m_{js}| \leq 1/s$, the error resulting from fixing the upper integration limit at $1000$ is $\int_{1000}^{\infty} r(s,x) ds < \int_{1000}^{\infty} 1/s^2 ds = 10^{-3}$. However, it is not always possible to determine the rate of decay of the magnitude of the imaginary moments with only numerical values. So it is difficult to come up with a universal rule for the integration range that is suitable for all moment sequences.} 

For the GP method, we propose a dynamic way to adjust the values of the upper limit and step size. We start at $F_{\rm{GP}(0.1, 1000)}$. If the cdf determined by the previous selections exceeds $[0,1]$, we decrease the step size to $1/3$ of the previous one and increase the upper limit to $3$ times the previous one, and then recalculate the cdf. We repeat this process until the range of the cdf is in $[0,1]$. In some cases, the updated step size and upper limit may cause memory problems. In the averaging process described in \Cref{sec:perf}, we discard those cases if we detect that the required memory is larger than the available memory.

\begin{figure}
	\centering
	\subfloat[$\upsilon = 1000$]{\includegraphics[width = 0.48\figwidth]{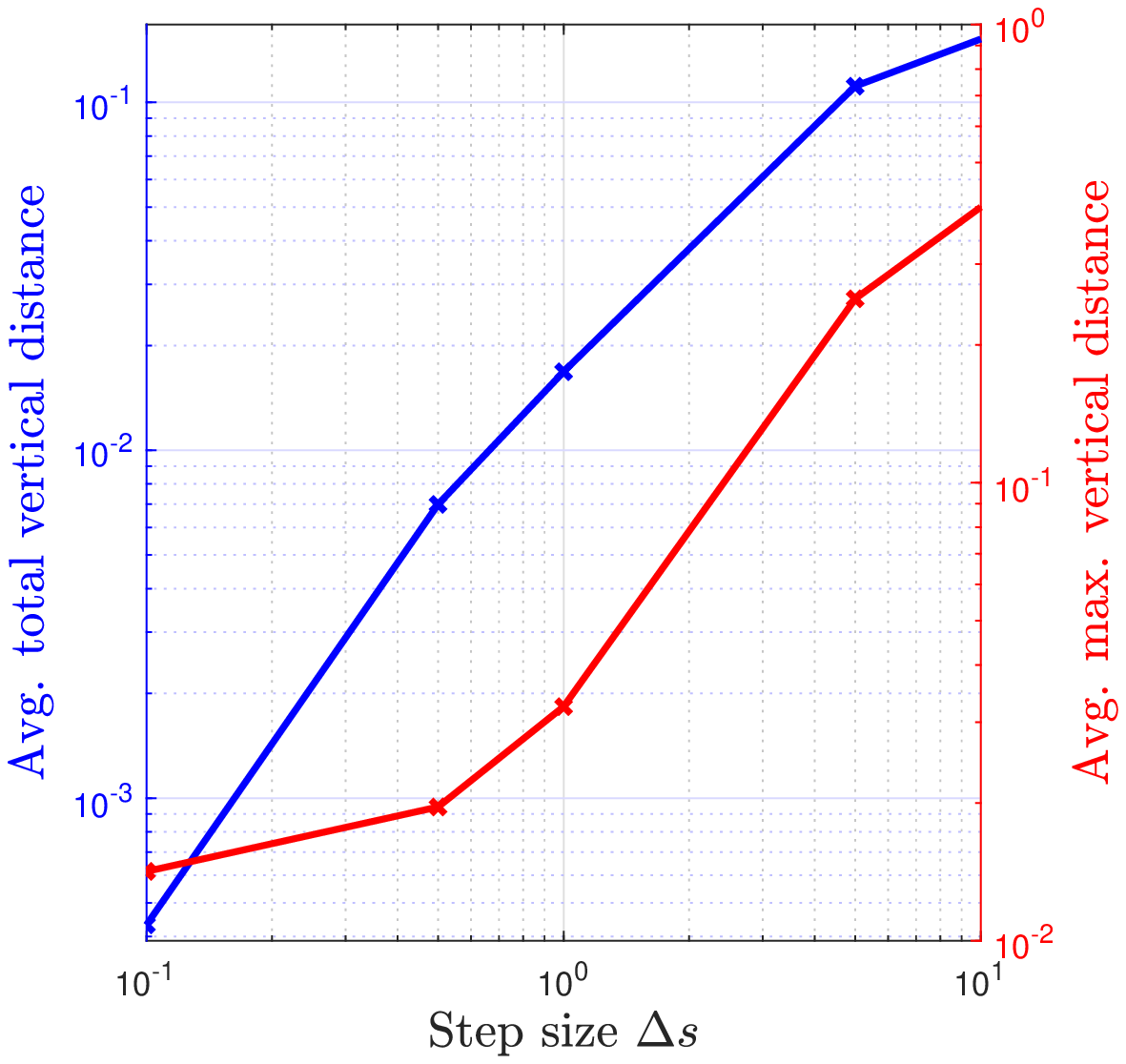}}\quad
	\subfloat[$\Delta s = 0.1$]{\includegraphics[width = 0.48\figwidth]{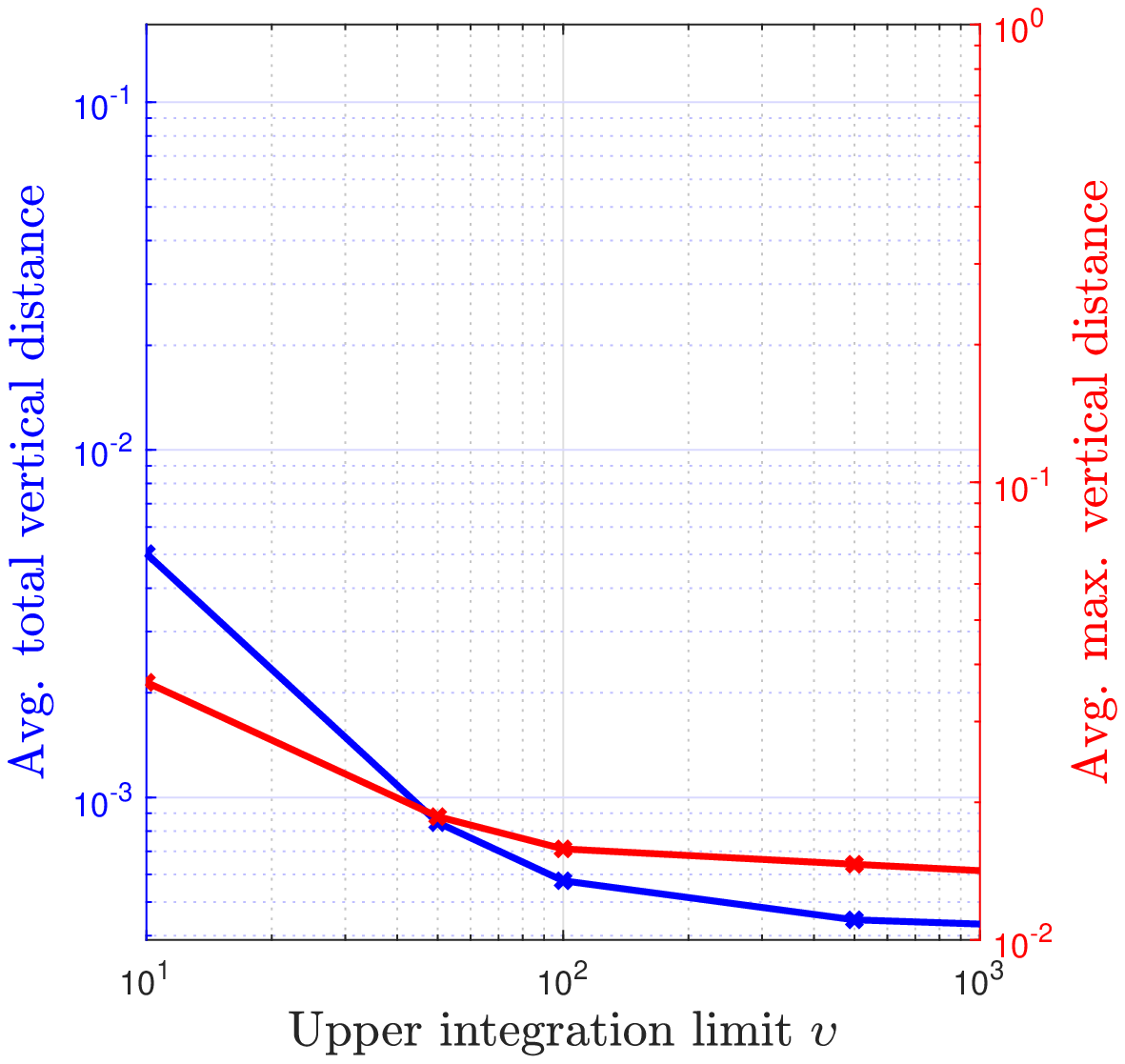}} 
	\caption{Average of the total and maximum distances between $F$ and polished $\hat{F}_{\rm{GP}(\Delta s,\upsilon )}$ which are averaged over $100$ randomly generated beta mixtures as described in \Cref{sec:bypdfs}. }\label{fig:gp_precision}
\end{figure}

\subsubsection{FJ method}

Equipped with a complete orthogonal system on $[0,1]$, denoted as $\left(T_k(x)\right)_{k = 0}^{\infty}$,  we can expand a function $l$ defined on $[0,1]$ to a generalized Fourier series as
\begin{align}
	l(x) = \sum_{k = 0}^{\infty} c_k T_k (x).
\end{align}
The Jacobi polynomials $P_n^{(\alpha, \beta)}$ are such a class of orthogonal polynomials defined on $[-1,1]$. 
They are orthogonal w.r.t.~the weight function $\omega^{(\alpha, \beta)} (x) \triangleq  (1-x)^\alpha (x+1)^{\beta}$, i.e., 
\begin{align}
	\int_{-1}^{1}  \omega^{(\alpha, \beta)} (x) P_n^{(\alpha, \beta)} (x)P_m^{(\alpha, \beta)} (x) \,dx  
	= \frac{2^{\alpha + \beta +1}  
		\left(n+1\right)^{(\beta)}
		\delta_{mn}
	}{(2n+\alpha + \beta +1) 
		\left(n+\alpha+1\right)^{(\beta)}} , 
\end{align}
where $\delta_{mn}$ is the Kronecker delta function, and $a^{(b)}\triangleq{\frac {\Gamma (a+b)}{\Gamma (a)}}$ is the Pochhammer function.
For $\alpha = \beta = 0$, they are the Legendre polynomials; for $\alpha = \beta = \pm 1/2$, they are the Chebyshev polynomials. The shifted Jacobi polynomials
\begin{align}
	R_n^{(\alpha, \beta)} (x) &\triangleq P_n^{(\alpha, \beta)} (2x-1) = \sum_{j = 0}^{n} \binom{n + \alpha}{j} \binom{n+\beta}{n - j} (x-1)^{n-j} x^j,
\end{align}
are defined on $[0,1]$. Letting $\eta_n^{(\alpha, \beta)}\triangleq \frac{ \left(n+1\right)^{(\beta)} 
}{(2n+\alpha + \beta +1) 
	\left(n+\alpha+1\right)^{(\beta)}
}$, the corresponding orthogonality condition w.r.t.~the weight function $w^{(\alpha, \beta)}(x) \triangleq (1-x)^{\alpha} x^\beta$ is
\begin{align}
	\int_{0}^{1}  w^{(\alpha, \beta)}(x)  R_n^{(\alpha, \beta)} (x)R_m^{(\alpha, \beta)} (x) \,dx=  \eta_n^{(\alpha, \beta)}\delta_{mn}.
\end{align} 


When applying the FJ series to solve \eqref{eq:truncated}, the usual approach is to expand the target function (cdf or pdf) as 
\begin{align}
	w^{(\alpha, \beta)}(x) \sum_{k = 0}^{\infty} c_k R_k^{(\alpha, \beta)} (x),
\end{align}
so that the moments can be utilized in the orthogonality conditions. {However, an approximation obtained by the $n$-th partial sum of the expansion is usually not a solution to the THMP, for two reasons. First, the remaining Fourier coefficients $c_m, m>n$, may not all be zero. Second, the approximations may not preserve the properties of probability distributions. }

Convergence is a major concern when approximating a function by an FJ series. While no conditions are known that are both sufficient and necessary, there exist sufficient conditions and there exist necessary conditions, both of which depend on the values of $\alpha, \beta$ and the properties of the approximated functions \cite{pollard1972convergence,li1995pointwise}. For the THMP, a basic question is whether $\alpha, \beta$ should depend on the moment sequence or not and whether to approximate pdfs or cdfs. In the following, we will discuss in detail three different ways that have been proposed in \cite{shohat1950problem, guruacharya2018approximation, schoenberg1973remark}. In particular, the authors in  \cite{guruacharya2018approximation} set $\alpha, \beta$ depending on the moment sequences and only approximate pdfs.  In \cite{shohat1950problem, schoenberg1973remark}, cdfs and pdfs are approximated, respectively, and both fix $\alpha=\beta = 0$.
 
The FJ method based on the FJ series is proposed in \cite{guruacharya2018approximation}. $\alpha, \beta$ are chosen such that the first two Fourier coefficients $c_1$ and $c_2$ are zero, i.e., $\alpha+ 1 =  (m_1- m_2) (1-m_1)/{(m_2 - m_1^2)}$ and $\beta + 1 =  {(\alpha + 1) m_1}/{(1 - m_1)}$. For any positive integer $n$, the approximation by the FJ method is defined as follows. 
\begin{Def}[Approximation by the FJ method \cite{guruacharya2018approximation}]
\begin{align}
	f_{\mathrm{FJ},n}(x)  \triangleq w^{(\alpha, \beta)} (x) \sum_{ k = 0}^n c_k R_k^{(\alpha, \beta)} (x), 
\end{align} 
where, by the orthogonality condition, 
\begin{align}
	c_m  
	= & \frac{1}{\eta_m^{(\alpha, \beta)}}  \sum_{ j = 0}^m \binom{m+\alpha}{j} \binom{m+\beta}{m-j} \sum_{ k = 0}^{m-j} \binom{m-j}{k} (-1)^k m_{m-k}.
\end{align}  
The corresponding cdf is 
\begin{align}
	F_{\mathrm{FJ},n}(x)  
	= \frac{
		\left(\alpha+1\right)^{(\beta+1)}
	}{ 
		\Gamma( \beta+1)  }  \int_{0}^{x} w^{(\alpha, \beta)} (t) dt  
	- \sum_{ k = 1}^n \frac{c_k}{k} w^{(\alpha+1, \beta+1)}(x)  R_{k-1}^{(\alpha+1, \beta+1)} (x). 
\end{align}
\end{Def} 

  $F_{\mathrm{FJ},1} = F_{\mathrm{FJ},2}$ correspond to the beta approximation \cite{haenggi2016meta}, where $m_1$ and $m_2$ are matched with the corresponding moments of a beta distribution.

There are three drawbacks of the FJ method. First, as $\alpha, \beta$ vary significantly for different moment sequences, the convergence properties of the FJ method vary strongly as well. {So it is impossible to give a conclusive answer to whether such approximations converge or not.} It is a critical shortcoming of the FJ method that it only converges in some cases. Second, as $\alpha$ and $\beta$ are generally not zero, the approximate pdf is always $0$ or $\infty$ at $x = 0$ and $x = 1$, which holds only for a limited class of pdfs. Last but not least, it is possible that, as a functional approximation, the FJ method leads to non-monotonicity in the approximate cdf. {We address this by applying the tweaking mapping in \Cref{def:correction} to $F_{\mathrm{FJ},n} |_{\mathcal{ U}_n}$.}

\subsubsection{FL method} \label{sec:FL}
In \cite{shohat1950problem, schoenberg1973remark}, the authors suggest fixing $\alpha = \beta = 0$ in advance. The resulting polynomials are Legendre polynomials with weight function $1$. To expand a function $l$ by the FL series, a necessary condition for convergence is that the integral $\int_0^1 l(x) (1-x)^{-1/4} x^{-1/4} \, dx$ exists \cite{szeg1939orthogonal}. Furthermore, \cite{pollard1972convergence} shows that if $l \in \mathcal{L}^p$ with $p > 4/3$, the FL series converges pointwise. There still remains the question whether we should approximate the pdf or the cdf. In the following, we will discuss the two cases. {The first one is considered in \cite{schoenberg1973remark}, and the second one is mentioned in \cite{shohat1950problem} but the concrete expressions of the coefficients are missing. Though the two approaches have been proposed, to the best of our knowledge, their convergence properties have not been compared.}

We denote the corresponding pdf of $F$ as $f$ and $R_k^{(0,0)}$ as $R_k$.  

\paragraph{\emph{Approximating the pdf \cite{schoenberg1973remark}} }  
The $n$-th partial sum of the FL expansion of the pdf is 
\begin{align}
	\sum_{k = 0}^{n} c_k R_k(x),
\end{align}
where, by the orthogonality condition,  
\begin{align}
	c_m 
	&= (2m+1) \sum_{j = 0}^m \binom{m}{j} \binom{m}{m - j} \sum_{k =0}^{m-j} \binom{m-j}{k}(-1)^{k} m_{m-k}.
\end{align}

It does not always hold that $f \in \mathcal{L}^p[0,1]$ with $p > 4/3$. Thus, the series does not converge for some $f$, for example, the pdf of $\mathrm{Beta} (1, \beta)$ with $0< \beta \leq 1/4$ or $\mathrm{Beta} (\alpha, 1)$ with $0< \alpha \leq 1/4$.  Worse yet, it is possible that the FL series diverges everywhere, e.g., the FL series of $f(x) =  (1-x)^{-3/4}/4$.

{So similar to the FJ method, approximating the pdf does not guarantee convergence. }

\paragraph{\emph{Approximating the cdf \cite{shohat1950problem}}} The $n$-th partial sum of the FL expansion of the cdf is
\begin{align}
	\sum_{k = 0}^{n} c_k R_k(x), 
\end{align} 
where, by the orthogonality condition, 
\begin{align}
	c_m 
	&= (2m+1) \sum_{j = 0}^m \binom{m}{j} \binom{m}{m - j}    \sum_{k =0}^{m-j} \binom{m-j}{k}(-1)^{k} 
	\frac{1 - m_{m-k+1}}{m-k+1}.\label{eq:flcm}
\end{align}
Since $F$ is always bounded and thus $F \in \mathcal{L}^p[0,1]$ for any $p \geq 1$, pointwise convergence holds \cite{pollard1972convergence}.

\begin{figure}
	\centering  
	{\includegraphics[width=0.6\figwidth]{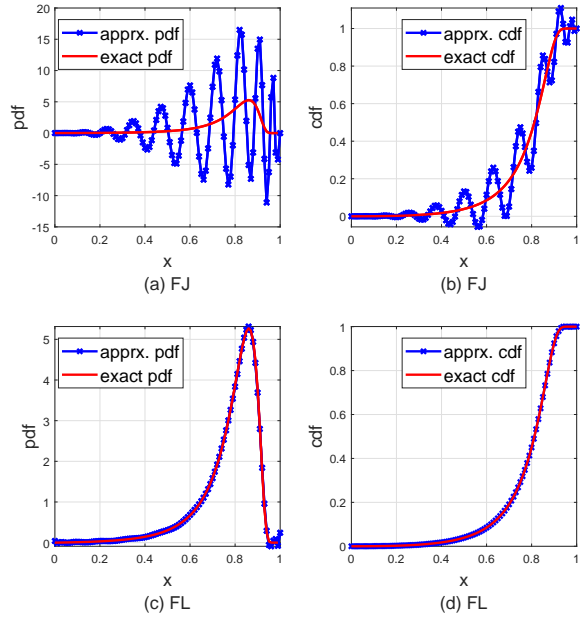}}
	\caption{ The upper two plots are for the 
		FJ method \cite{guruacharya2018approximation} with $n = 20$, and the lower two plots are for the FL series with $n = 20$. The cdf is given in \eqref{eq:accdf}. }\label{fig:FJcon20}
\end{figure}

\begin{figure}
	\centering  
	{\includegraphics[width=0.6\figwidth]{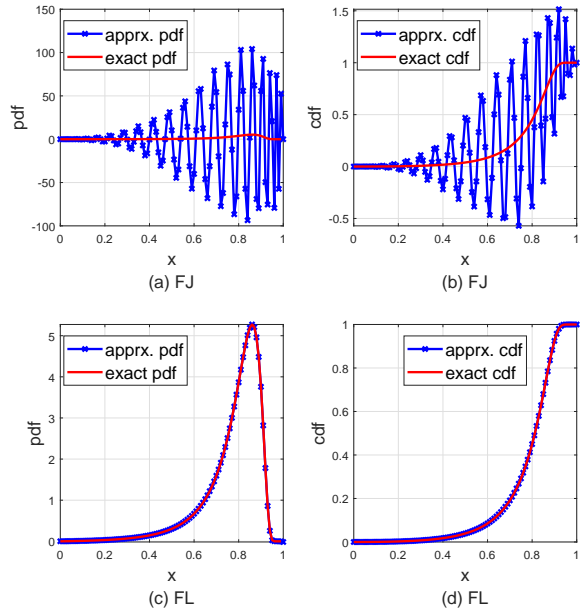}}
	\caption{  The upper two plots are for 
		FJ method \cite{guruacharya2018approximation} with $n = 40$, and the lower two plots are for the FL series with $n = 40$. The cdf is given in \eqref{eq:accdf}. }\label{fig:FJcon40}
\end{figure}

\paragraph{\emph{Comparisons of the FJ method and  approximating the pdf via the FL series}}
In \Cref{fig:FJcon20,fig:FJcon40}, where 
\begin{align}
	{F}(x) = 1- \frac{\exp\left(- \frac{x^2}{50(1-x)^2}\right)  }{  1 +  \frac{x^2}{50(1-x)^2} }\label{eq:accdf},
\end{align}
we approximate the pdf by the FJ method and by the FL series. We observe that in both cases,  the FL approximations are much better than the FJ ones. Also, the FL approximations do better as $n$ increases.

\paragraph{\emph{Comparisons of approximating the pdf and cdf by the FL series}}
{\Cref{fig:FLcon} shows the approximate cdfs by using the FL series to approximate the pdf and then integrate and to directly approximate the cdf, respectively.}
We observe that approximating cdfs leads to a more accurate result.

\begin{figure}
	\centering   
	\includegraphics[width = 0.45\figwidth]{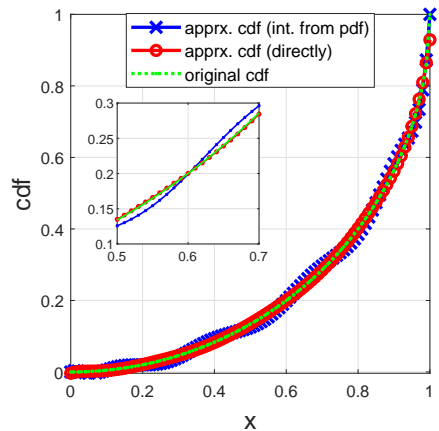}
	\caption{  Reconstructed cdfs of $F(x) = 1 - \sqrt{1-x^2}$ based on $10$ moments. The blue curve shows the reconstructed cdf from approximating the pdf, the red curve shows the reconstructed cdf from approximating the cdf, and the green curve shows $F$.  
	}\label{fig:FLcon}
\end{figure}

Considering that a pdf may not be in $\mathcal{L}^p[0,1]$ with $p > 4/3$ but a cdf always is, in the following, we use the FL series to approximate the cdf and call it the FL method. 
For any positive integer $n$, the approximation by the FL method is defined as follows. 
\begin{Def}[Approximation by the FL method]
\begin{align}
	F_{\mathrm{FL},n}(x) \triangleq& 
	\sum_{m = 0}^{n} c_m R_m(x),
\end{align}
where $c_m, \, m \in [n]_0$, are given by \eqref{eq:flcm}. 
\end{Def}
Since it is a functional approximation, we apply the tweaking mapping in \Cref{def:correction} to $F_{\mathrm{FL},n-1}|_{\mathcal{U}_{n}}$.

\subsection{New methods}
\subsubsection{FC method}
{Similar to the idea in \Cref{sec:FL} that we fix $\alpha$ and $\beta$ in advance, here we use the FC series for approximation, i.e., we set $\alpha = \beta = -1/2$.}  
To guarantee convergence, {similar to the FL method}, we also focus on approximating the cdf and call it the FC method. For any positive integer $n$, the approximation by the FC method is defined as follows. 

\begin{Def}[Approximation by the FC method]
\begin{align}
	F_{\mathrm{FC}, n} (x)\triangleq 
	\begin{cases}
		x^{-\frac{1}{2}}\left(1-x\right)^{-\frac{1}{2}}  \sum_{m = 0}^n c_m R_m^{\left(-\frac{1}{2}, -\frac{1}{2} \right)}(x), \quad x \in (0,1),\\
		x, \quad x \in \{0,1\},
	\end{cases}
\end{align}
where, by the orthogonality condition, 
\begin{align}
	c_m 
	&=  c'_m  \sum_{j = 0}^m \binom{m-1/2}{j} \binom{m-1/2}{m - j}    \sum_{k =0}^{m-j} \binom{m-j}{k}(-1)^{k} 
	\frac{1 - m_{m-k+1}}{m-k+1},\label{eq:fccm}
\end{align}
$c'_0 = 1/\pi$ and $c'_m = 2 \left( \frac{\Gamma(m+1)}{\Gamma(m+1/2)}\right)^2$, $m \in [n]$.
\end{Def}
As $F$ and $x^{\frac{1}{2}}\left(1-x\right)^{\frac{1}{2}} $ are of bounded variation, the function $F x^{\frac{1}{2}}\left(1-x\right)^{\frac{1}{2}}$ 
approximated by the FC series 
has bounded variation, thus uniform convergence holds \cite{mason2002chebyshev}. {To preserve the properties of probability distributions, we apply the tweaking mapping in \Cref{def:correction} to $F_{\mathrm{FC}, n-1}|_{\mathcal{ U}_{n}}$. }

%
%

\subsubsection{CM method}
As the CM inequalities provide the infima and suprema at the points of interest among the family of solutions given a sequence of length $n$, it is important to analyze the distribution of the unique solution given the sequence of an infinite length. Here, we explore this distribution at point $x = 0.5$ by simulating the distribution of infima and suprema of a sequence of length $8$ within the initial ranges define by the infima and suprema given the first three elements in the sequence. For example, let $\underaccent{\bar}{F}_3$ and  $\bar{F}_3$ denote the infimum and supremum defined by the CM inequalities at point $0.5$ given a sequence $(m_1, m_2, m_3)$, and let  $\underaccent{\bar}{F}_8$ and  $\bar{F}_8$  denote the infimum and supremum defined by the CM inequalities at point $0.5$ given a sequence $(m_1, m_2, m_3, m_4, m_5, m_6, m_7, m_8)$. We focus on the relative position of $\underaccent{\bar}{F}_8$ and  $\bar{F}_8$  within the range defined by $\underaccent{\bar}{F}_3$ and  $\bar{F}_3$, i.e., the normalized  $x = \left(\underaccent{\bar}{F}_8 - \underaccent{\bar}{F}_3 \right)/\left(\bar{F}_3- \underaccent{\bar}{F}_3\right)$ and $x = \left(\bar{F}_8 - \underaccent{\bar}{F}_3 \right)/\left(\bar{F}_3 - \underaccent{\bar}{F}_3\right)$. The average is calculated as $\left(\bar{F}_8+ \underaccent{\bar}{F}_8 \right)/2$. 

\begin{figure}
	\centering  
	{\includegraphics[width=\figwidth]{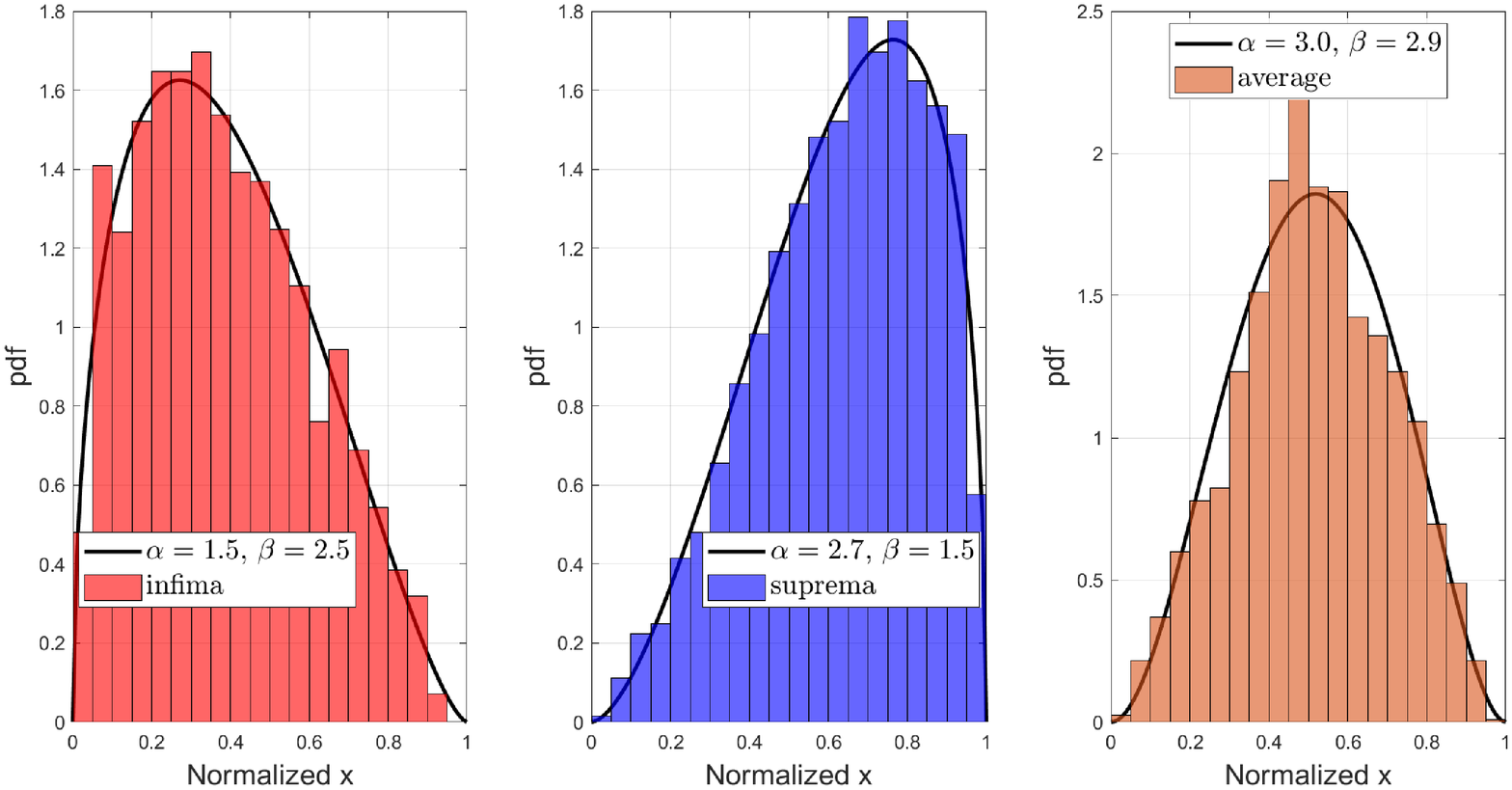}}
	\caption{The normalized histogram of infima, suprema, and average and the corresponding beta distribution fit for $50$ realizations of moment sequences of length $3$ generated by \Cref{sch:canonical}, each of which has been further extended $50$ times to  length $8$  by \Cref{sch:canonical}. The two parameters of the beta distribution fit are based on the first two empirical moments. }\label{fig:randomhist}
\end{figure}

\begin{figure}
	\centering  
	{\includegraphics[width=\figwidth]{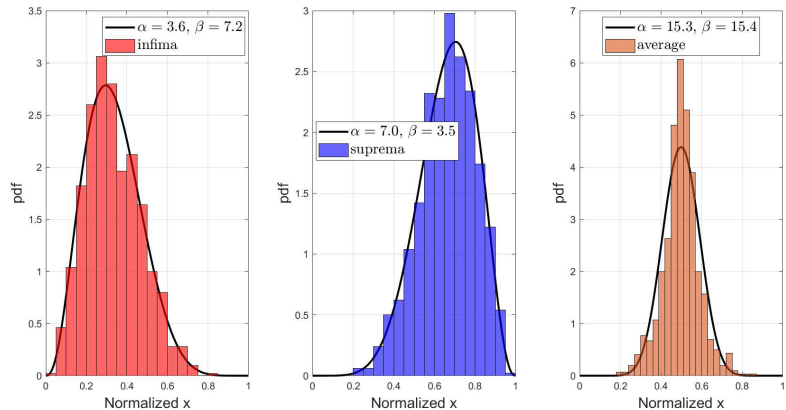}}
	\caption{The normalized histogram of infima, suprema, and average and the corresponding beta distribution fit for $1000$ realizations of moment sequences of length $8$ generated by \Cref{sch:canonical}, each of which has the same first three elements as $(1/2, 1/3, 1/4)$. The two parameters of the beta distribution fit are based on the first two empirical moments. }\label{fig:unifhist}
\end{figure}

\Cref{fig:randomhist} and \Cref{fig:unifhist} show the normalized histograms with beta distribution fit. From both of them, we can observe that there is symmetry between the distribution of the infima and suprema, and their average is concentrated in the middle. In particular, in \Cref{fig:unifhist}, we observe that its average concentrates more in the middle compared to \Cref{fig:randomhist}, since the first three moments come from the uniform distribution. This finding indicates that the unique solution of a moment sequence of infinite length is  likely to be centered within the infima and suprema defined by its truncated sequence.  Therefore, it is sensible to consider their average as an approximation to a solution to the THMP and we call it the CM method. For any positive integer $n$, the approximation by the CM method is defined as follows. 
\begin{Def}[Approximation by the CM method]
\begin{align}
	F_{\rm{CM},n} \triangleq \frac{1}{2} \left(\sup_{F\in \mathcal{F}_n} F(x) + \inf_{F\in \mathcal{F}_n} F(x) \right), 
\end{align}
where $\sup_{F\in \mathcal{F}_n} F(x) $ and $\inf_{F\in \mathcal{F}_n} F(x)$ are given by the CM inequalities. 
\end{Def}
From \Cref{fig:markovAvg}, we observe that the average is very close to the original cdf, and it is slightly better than the beta approximation. It is worth noting that by averaging, the error can be bounded: For all $x \in [0,1]$ and for all $F \in \mathcal{F}_n$, 
\begin{align}
	|F_{\rm{CM},n}(x) - F(x) | \leq \frac{1}{2} \left(\sup_{F\in \mathcal{F}_n} F(x) - \inf_{F\in \mathcal{F}_n} F(x) \right).
\end{align}

{Although $F_{\rm{CM},n}$ lies in the band defined by the infima and suprema from the CM inequalities, it is not necessarily a solution (among the infinitely many solutions) to the THMP. Moreover, it is almost never a solution. For example,  $F_{\rm{CM},1}$ is a solution only when $m_1 = 1/2$. 
}

\begin{figure}[!htbp]
	\centering
	\subfloat[$F(x) =\exp(-x/(1-x))$. The total difference is $0.015$ between $\hat{F}_{\rm{CM},6}$ and the actual cdf and $0.019$ for the beta approximation.]{\includegraphics[width =0.45 \figwidth]{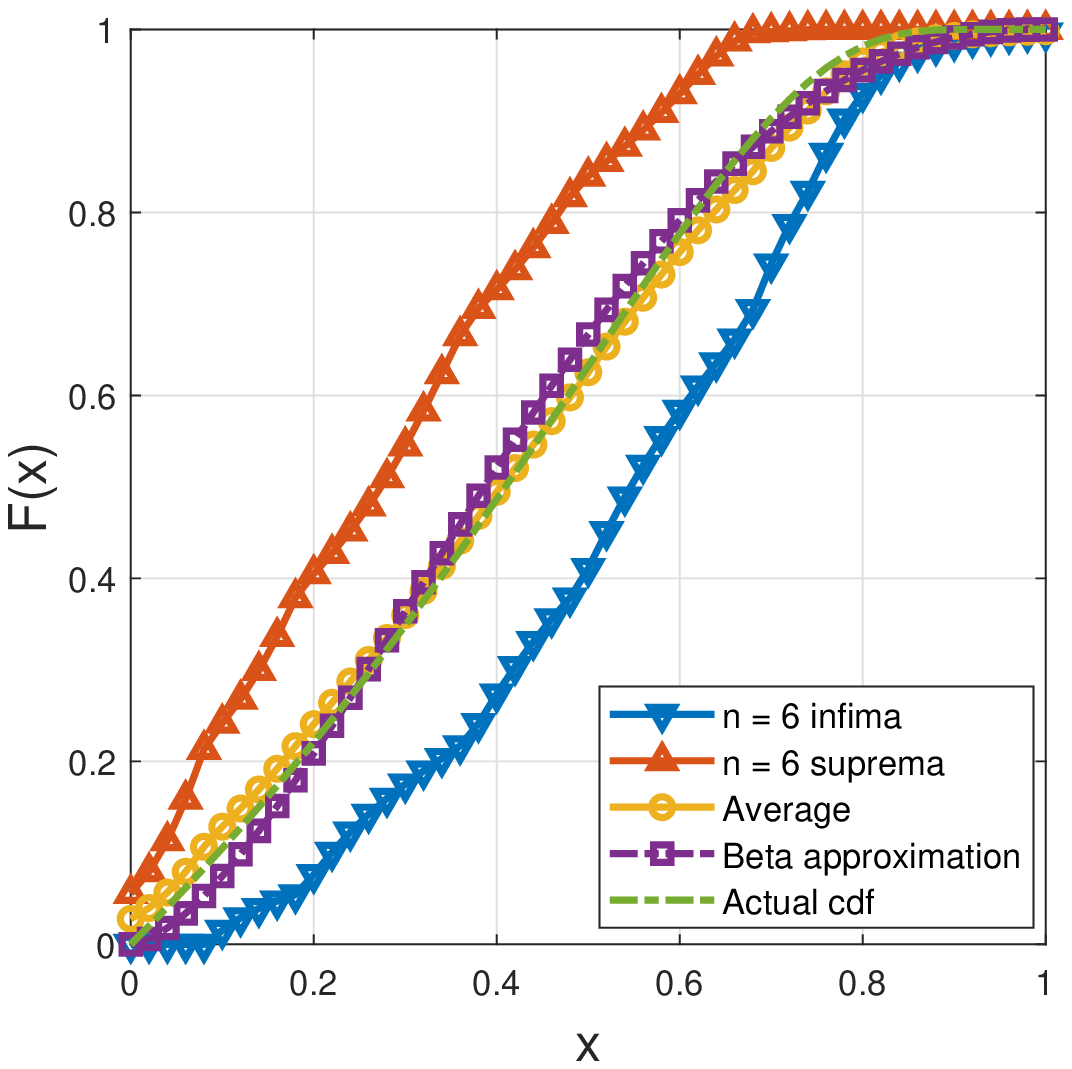}} \quad
	\subfloat[$F(x) = \frac{1}{2}\int_0^1 (x^{10}(1-x) + x(1-x)^{10}) \, dx$. The total difference is $0.0199$ between $\hat{F}_{\rm{CM},6}$ and the actual cdf and $0.0451$ for the beta approximation.]{\includegraphics[width =0.45 \figwidth]{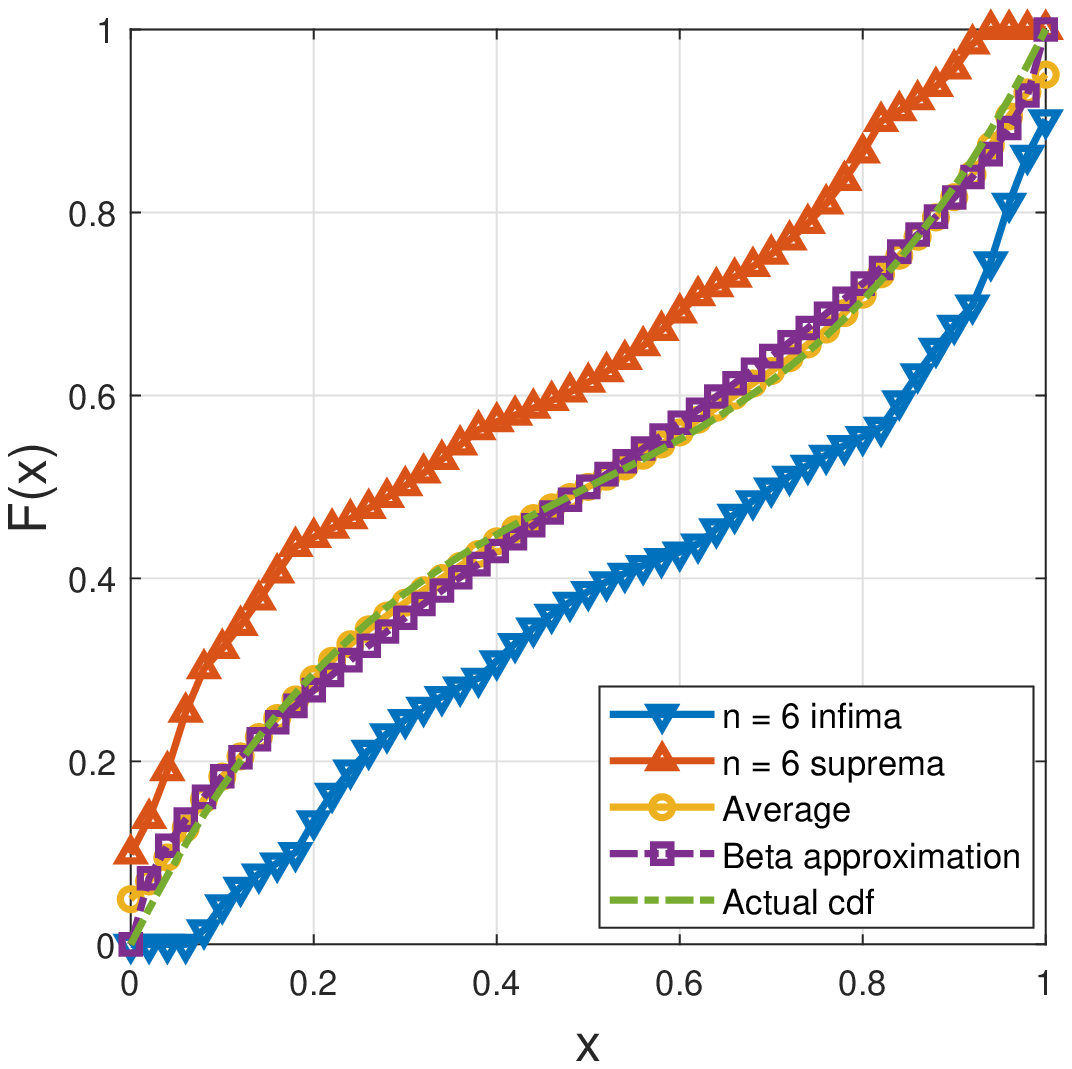}} 
	\caption{The infima and suprema from the CM inequalities, their average, and the beta approximation given $(m_k)_{k = 0}^6$. }\label{fig:markovAvg}
\end{figure}

\subsection{The methods as linear transforms}
The calculation of the coefficients in the BM method and methods based on the FJ series, such as the FL and FC method, where the coefficients $\alpha$ and $\beta$ are fixed, can be written as linear transforms of the given sequence and the transform matrix is fixed for each  method. This can be done offline and accelerates the calculation.\footnote{The FJ method can also be written as linear transform but it is meaningless to do so due to the undetermined coefficients $\alpha$ and $\beta$.} We discuss the BM and FL method in detail. 

As we recall, for the BM method \cite{haenggi2018efficient}, we can write $\mathbf{h} \triangleq  \left(h_l\right)_{l = 0}^n$ as a linear transform of $\mathbf{m} \triangleq  (m_l)_{l = 0}^{n}$, i.e., 
\begin{align}
	\mathbf{h} = \mathbf{A} \mathbf{m}, 
\end{align}
where $\mathbf{m}$ and $\mathbf{h}$ are column vectors and the transform matrix $\mathbf{A} \in \mathbb{Z}^{(n+1)\times(n+1)}$ is given by 
\begin{align}
	A_{ij} \triangleq \binom{n}{j} \binom{j}{i} (-1)^{j-i} \mathbbm{1}{(j \geq i)}, \quad i,j \in [n]_0, 
\end{align}
where $\mathbbm{1}$ is the indicator function. The transform matrix needs to be calculated only once for each $n$, which can be done offline. The matrix-vector multiplication requires $\frac{(n+1)(n+2)}{2}$ multiplications. 

Similarly, for the FL method, by reorganizing \eqref{eq:flcm}, we can write $\mathbf{c} \triangleq  \left(c_l\right)_{l = 0}^n$ as a linear transformation of $\mathbf{\hat{m}} \triangleq  (m_l)_{l = 1}^{n+1}$, i.e., 
\begin{align}
	\mathbf{c}= \mathbf{\hat{A}} (\mathbf{1} - \mathbf{\hat{m}}), \label{eq:fllt}
\end{align}
where $\mathbf{c}$ and $\mathbf{\hat{m}}$ are understood as column vectors, $\mathbf{1}$ is the 1-vector of size $n+1$, and the transform matrix $\mathbf{\hat{A}} \in \mathbb{Z}^{(n+1)\times(n+1)}$ is given by 
\begin{align}
	\hat{A}_{k l} \triangleq \frac{(-1)^{(k-l)}(2k+1)}{l+1} \sum_{j = 0}^{l} \binom{k}{j} \binom{k}{k-j} \binom{k-j}{k-l}\mathbbm{1}(l\leq k), \quad k, l \in [n]_0.
\end{align}
$\mathbf{\hat{A}}$ is a lower triangular matrix. As for the BM method, the transform matrix needs to be calculated only once for the desired level of accuracy, and the above matrix-vector multiplication requires $\frac{(n+1)(n+2)}{2}$ multiplications. More details about the accuracy requirements will be given in \Cref{sec:accuracyHMT}.

\section{Generation of moment sequences and performance evaluation}
To compare the performance of different HMTs, moment sequences and corresponding functions serving as the ground truth are needed. The problem of generating representative moment sequences is non-trivial since the probability that sequences i.i.d.~ on $[0,1]^n$ are moment sequences is about $2^{-n^2}$ \cite[Eq.~1.2]{chang1993normal}. In this section, we provide  three different ways to solve it: The first one generates uniformly distributed integer moment sequences; the other two generate moments from either analytic expressions of moments or analytic expressions of cdfs. In these two cases, the parameters are randomized. We use the exact cdf as the reference function if it is given by analytic expressions. Otherwise, we use the GP approximations as the reference functions.

\subsection{Moment sequences and their rates of decay}\label{sec:rateofdecay}

As discussed in \Cref{sec:cm}, the CM method is able to fully recover discrete distributions with a finite number of jumps, while other methods, such as the ME method, are unable to do that. This indicates that different HMTs may have different performance for different types of distributions, classified by the decay order of the sequence. We first give some examples of moment sequences of common distributions.
\begin{Exam}[Beta distributions]\label{exam:beta} 
	For $X \sim \operatorname{Beta}(\alpha, \beta)$, 
	\begin{align}
		m_n & =\prod_{r = 0}^{n-1} \frac{\alpha +r}{\alpha + \beta +r}  = \frac{\Gamma(\alpha + n ) }{\Gamma(\alpha ) } \frac{\Gamma(\alpha + \beta ) }{\Gamma(\alpha + \beta +n  ) }. 
	\end{align} 
	By Stirling's formula, as $n \to \infty$,  
	\begin{align}
		m_n & \sim e^\beta\frac{\Gamma(\alpha + \beta ) }{\Gamma(\alpha ) }\sqrt{\frac{\alpha + n-1}{\alpha + \beta + n-1}} \frac{(\alpha + n-1)^{\alpha + n-1}  }{  (\alpha + \beta + n-1)^{\alpha + \beta + n-1}}, \\
		& \sim e^\beta \frac{\Gamma(\alpha + \beta  ) }{\Gamma(\alpha ) } \frac{(\alpha + n-1)^{\alpha + n-1}}{(\alpha + \beta + n-1)^{\alpha + \beta + n-1}}, \\
		& \sim   n^{-\beta}, \quad \text{as } \left(1+\frac{\beta}{x}\right)^{x} \to e^\beta, \, x\to \infty. 
	\end{align}
	Thus, $m_n = \Theta(n^{-\beta})$.
\end{Exam}
 
\begin{Exam}[]
	Let the density function be $f(x) = 2 x/s^2  \mathbf{1}_{[0,s]}(x)$ with $0<s< 1$. Then $m_n =  2 s^n /(n+2) = O(\exp(-\log(1/s) n )) = \Omega(\exp(-s_1 n ))$ for all $s_1 > \log(1/s)$.
\end{Exam}

Given the examples shown above and that the performance of different HMTs is likely to be dependent on the rate of decay (or type of decay) of the moment sequences, we would like to develop a method to categorize different moment sequences based on their tail behavior. To this end, we introduce six types of moment sequences depending on their rates of decay.

\begin{Def}[Rate of decay]\label{def:rateofdecay}
	\leavevmode
	\begin{enumerate} 
		\item[1)] Sub-power-law decay: In this case, for all $s >0$, $m_n = \Omega(n^{-s})$.
		\item[2)] Power-law decay: In this case, there exists $s > 0$ such that $m_n = \Theta(n^{-s})$.
		\item[3)] Soft power-law decay: In this case, $m_n \neq \Theta(n^{-s}), \,\forall s \geq 0$ but there exist $s_1> s_2 \geq 0$ such that $m_n = \Omega(n^{-s_1})$ and $m_n = O(n^{-s_2})$.
		\item[4)]  Intermediate decay: In this case, for all $s_1 \geq 0$, $m_n = O(n^{-s_1})$ and for all $s_2 > 0$, $m_n = \Omega(\exp(-s_2 n))$.
		\item [5)] Soft exponential decay: In this case, $m_n \neq \Theta(\exp(-s n)), \, \forall s > 0$ but there exists $s_1> s_2>0$, $m_n = \Omega(\exp(-s_1 n))$ and $m_n = O(\exp(-s_2 n))$.
		\item[6)] Exponential decay: In this case, there exists $s>0$ such that $m_n = \Theta(\exp(-s n))$. 
	\end{enumerate}
\end{Def}
\subsection{Uniformly random (integer) moment sequences}
\subsubsection{Generation}\label{sec:cano}

Let $\Lambda$ denote the set of probability measures on $[0,1]$. For $\lambda \in \Lambda$ we denote by
\begin{align}
	m_n(\lambda) \triangleq \int_{0}^{1} x^n \, d\lambda, \quad n \in \mathbb{N}_0,
\end{align}
its $n$-th moment. Trivially, $m_0(\lambda) = 1$, $\forall \lambda \in \Lambda$. The set of all moment sequences is defined as
\begin{align}
	\mathcal{M} \triangleq \{\left(m_1(\lambda), m_2(\lambda), ...\right), \lambda \in \Lambda\}.
\end{align}
The set of the projection onto the first $n$ coordinates of all elements in $\mathcal{M}$ is denoted by $\mathcal{M}_n$. 
 
Now we recall how to establish a one-to-one mapping from the moment space $\interior \mathcal{M}$ onto the set $(0,1)^\mathbb{N}$ \cite{chang1993normal}, where $\interior \mathcal{M}$ denotes the interior of the set $\mathcal{M}$.
For a given sequence $\left(m_k\right)_{k \in [n-1]} \in \mathcal{M}_{n-1}$, the achievable lower and upper bounds for $m_n$ are defined as 
\begin{align}
	m_n^- \triangleq \min \{ m_n(\lambda), \lambda \in \Lambda \text{ with } m_k(\lambda) = m_k \text{ for } k \in [n-1]\}, \\
	m_n^+ \triangleq \max \{ m_n(\lambda), \lambda \in \Lambda \text{ with } m_k(\lambda) = m_k \text{ for } k \in [n-1]\}.
\end{align}
Note that $m_1^-  =0 $ and $m_1^+ = 1$. $m_n^-$ and $m_n^+$, $n \geq 2$, can be obtained by taking the infimum and supremum of the range obtained from solving the two linear inequalities $0\leq \underline{H}_{n}$ and $  0\leq \overline{H}_{n}$.

We now consider the case where $\left(m_k\right)_{k \in [n-1]}  \in \interior \mathcal{M}_{n-1}$, i.e., $m_n^+ >  m_n^-$.
The canonical moments are defined as 
\begin{align}
	p_n \triangleq \frac{m_n - m_n^-}{m_n^+ - m_n^-}. \label{eq:def}
\end{align}
By \cite{dette1997theory}, $\left(m_k\right)_{k \in [n]} \in \interior \mathcal{M}_n$ if and only if $p_n \in (0,1)$. In this way, $p_n$ represents the relative position of $m_n$ among all the possible $n$-th moments given the sequence $\left(m_k\right)_{k \in [n-1]} $.  
The following lemma \cite{chang1993normal} shows the relationship between the distribution on the moment space $\mathcal{M}_n$ and the distribution of $p_k$, $k \in [n]$. 

\begin{Lem}[{\cite[Thm.~1.3]{chang1993normal}}]\label{lem:chuang}
	For a uniformly distributed random vector $\left(m_k\right)_{k \in [n]} $ on the moment space $\mathcal{M}_n$, the corresponding canonical moments $p_k$, $k \in [n]$, are independent and beta distributed as
	\begin{align*}
		p_k \sim \mathrm{Beta}(n - k +1, n-k+1), \quad k \in [n].
	\end{align*}
\end{Lem}
Since $\mathcal{M}_n$ is a closed convex set, its boundary has Lebesgue measure zero. Therefore, although \eqref{eq:def} only defines the mapping on $\interior \mathcal{M}_n$, the $p_k$ are still well defined.

\Cref{lem:chuang} provides a method to generate a uniformly random moment sequence. The details are summarized in the following algorithm.

\begin{algorithm}[H]
	\caption{Algorithm for random moment sequences by canonical moments}\label{sch:canonical}
	\begin{algorithmic}[1]
		\renewcommand{\algorithmicrequire}{\textbf{Input:} }
		\renewcommand{\algorithmicensure}{\textbf{Output:} }
		\REQUIRE The number of moments $N$
		\ENSURE  $(m_n)_{n \in [N]}$ 
		\STATE Set $m_1^- = 0$, $m_1^+ = 1$ and $k = 1$.
		\WHILE {$j \leq N$}
		\STATE $p_k \sim \mathrm{Beta}(n-k+1, n-k+1)$, $m_k = p_k(m_k^+ - m_k^-) + m_k^-$ and $k = k+1$.  
		\STATE Calculate $m_k^-$ and $m_k^+$ based on $0\leq \underline{H}_{k}$ and $  0\leq \overline{H}_{k}$.
		\ENDWHILE 
	\end{algorithmic}
\end{algorithm}

\begin{figure}[!htbp]
	\centering
	\subfloat[Raw curves.\label{fig:recona}]{\includegraphics[width =0.48 \figwidth]{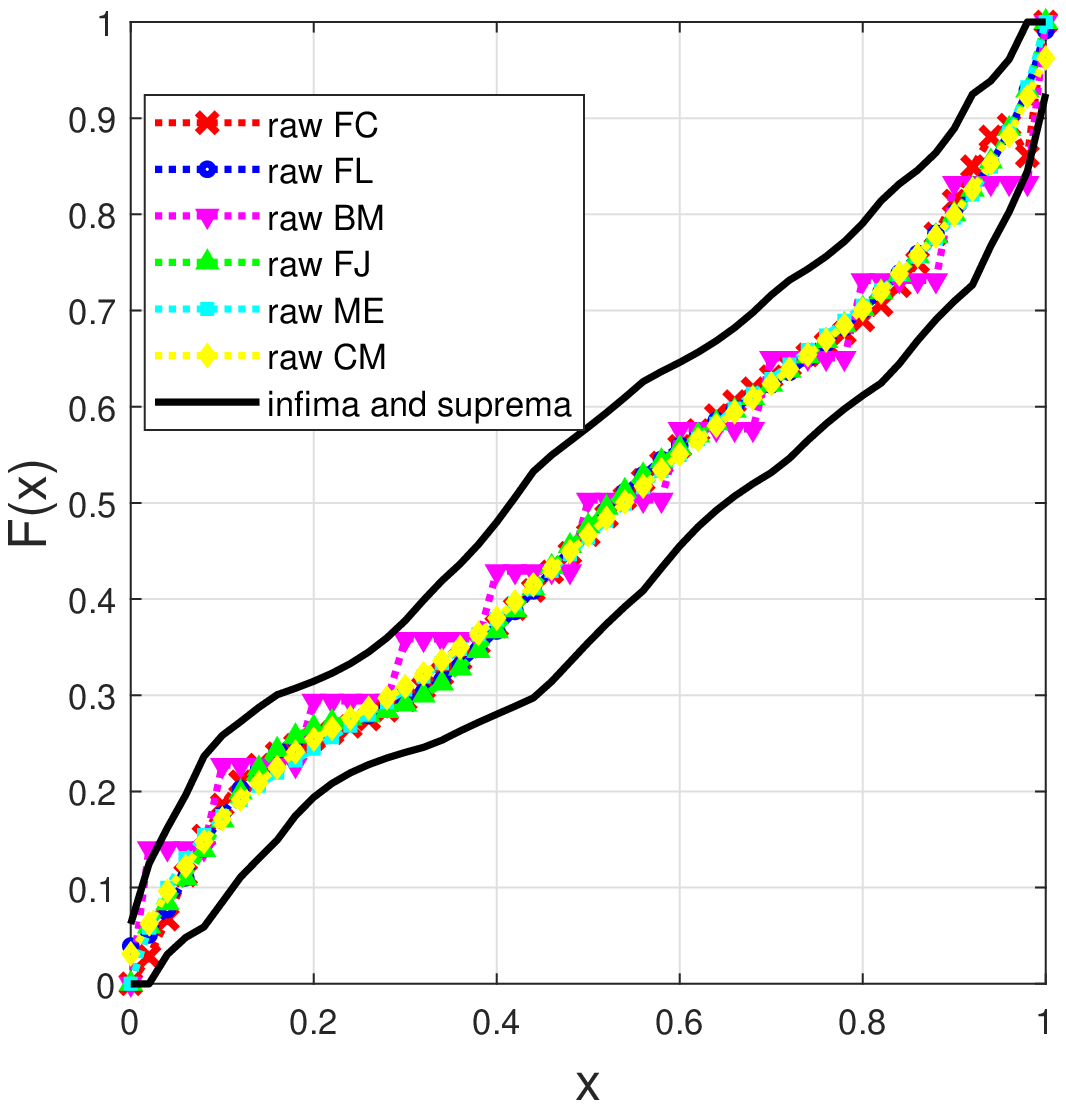}}\quad
	\subfloat[Polished curves.\label{fig:reconb}]{\includegraphics[width =0.48 \figwidth]{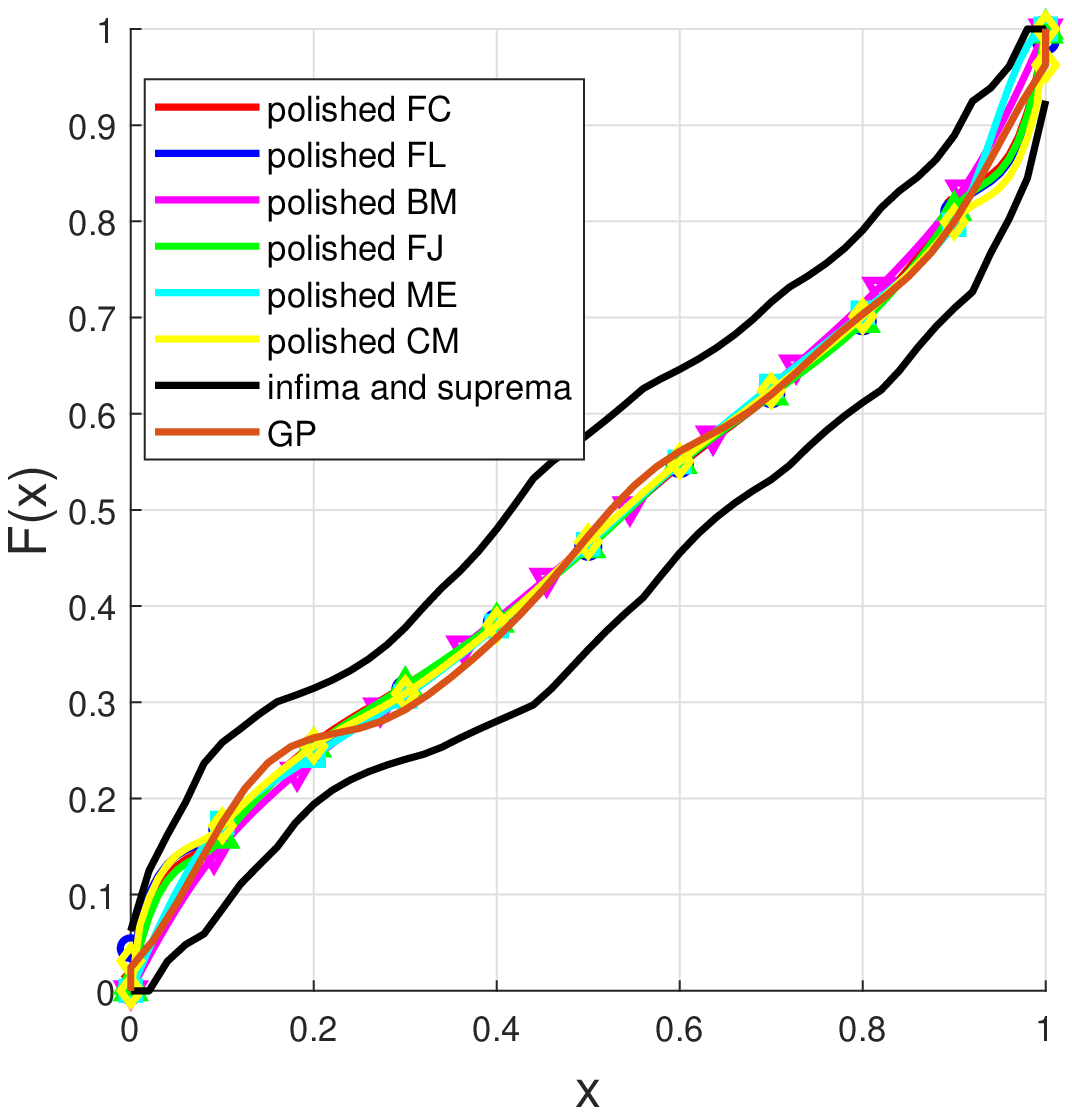}} 
	\caption{An example of the reconstructed cdf of each method. The $10$ moments are randomly generated by the canonical moments and the moments are $0.5256$, $0.3893$, $0.3188$, $0.2747$, $0.2442$, $0.2215$, $0.2038$, $0.1895$, $0.1775$ and $0.1674$. The polished curves are obtained after the tweaking mapping in \Cref{def:correction} as well as the monotone cubic interpolation for $\hat{F}_{\rm{BM},10}|_{\mathcal{ U}_{11}}$, $\hat{F}_{\rm{ME},10}|_{\mathcal{ U}_{10}}$, $\hat{F}_{\rm{FJ},10}|_{\mathcal{ U}_{10}}$, $\hat{F}_{\rm{FL},{9}}|_{\mathcal{ U}_{10}}$, $\hat{F}_{\rm{FC},9}|_{\mathcal{ U}_{10}}$, and $\hat{F}_{\rm{CM},10}|_{\mathcal{ U}_{10}}$.} 
\end{figure}

\subsubsection{Usage}
\Cref{fig:recona} shows an example of the reconstructed cdfs (only raw versions) of each method where moments are generated by \Cref{sch:canonical} and 10 moments are used. But since a reference function serving as the ground truth is missing, we cannot compare the performance of each method.
 It is natural to look for a way to reconstruct a solution to the HMT via (infinitely many) integer moments. However, no known reconstruction methods exists for integer moments. To address this problem, we take a detour to approximate the characteristic function first, and then apply the Gil-Pelaez theorem to reconstruct the cdf \cite{gil1951note}. Let $X \in [0,1]$ denote the random variable and $(m_n)_{n=0}^{\infty}$ denote its moments. The characteristic function of $X$ is
\begin{align}
	\phi_X(s)  & =  \mathbb{E}[\exp(jsX)], \quad s\in \mathbb{R}^+.
\end{align} 
Written in terms of its Taylor expansion,
\begin{align}
	\phi_X(s) & = \mathbb{E}\left[\sum_{k = 0}^{\infty} \frac{(jsX)^k}{k!}\right].
\end{align} 
By the dominated convergence theorem, we can change the order of expectation and summation, i.e.,  for any $s \in [0,\infty)$
\begin{align}
	\mathbb{E}\left[\sum_{k = 0}^{\infty} \frac{(jsX)^k}{k!}\right]   & = \sum_{k = 0}^{\infty} \frac{j^k s^k m_k}{k!}.\label{eq:dct}
\end{align} 
 
	The imaginary moments of $X$, $m_{js}$, are the characteristic function of $\log(X)$. When applied to $X \sim \exp(1)$, the expansion in \eqref{eq:dct} yields the gamma function, which is consistent with the analytic continuation property of the moments.

The truncated sum
\begin{align}
	\hat{\phi}_X(s)  = \sum_{k = 0}^{N} \frac{j^k s^k m_k}{k!} 
\end{align} 
is an approximation of the characteristic function $\phi_X(s)$, and, in turn, the integral 
\begin{align}
	\hat{F}_X(x) & \triangleq \frac{1}{2} - \frac{1}{\pi} \int_0^\infty  \frac{\Im{(e^{-jsx}\hat{\phi}_X(s))}}{s} ds \label{eq:gp_charac}, \quad x \in [0,1],
\end{align}
is an approximation of the exact cdf
\begin{align}
	{F}_X(x) & = \frac{1}{2} - \frac{1}{\pi} \int_0^\infty  \frac{\Im{(e^{-jsx}{\phi}_X(s))}}{s} ds \label{eq:gp_charac}, \quad x \in [0,1].
\end{align}
Therefore, by evaluating (approximating)  $\hat{F}_X$, we obtain an approximation of the exact cdf ${F}_X$. 

\Cref{fig:reconb} shows the reconstructions by the above method and all the HMTs. We observe that all the HMTs approximate well. However, even though the series converges for $s \in [0,\infty)$, the evaluation of each element with alternating signs could cause serious numerical problems, as $s$ increases. For each $s$, we need at least $N = \ceil{se}$ to approximate  $\phi_X(s)$ so that the absolute error is within $\sqrt{2/(\pi N)}$. Meanwhile, an accurate evaluation of $\hat{F}_X$ usually requires the calculation of the integrand up to $s = 1000$. This requires the calculation of $N = 2719$ moments, which is time-consuming and not practical for the performance evaluation purposes in \Cref{sec:perf}.

\subsection{Generating moments by analytic expressions} \label{sec:momentsbyanalytic}
As we need to generate moment sequences of variable lengths, it is natural to directly start with the infinite sequences $(m_n)_{n=0}^{\infty}$ and truncate them to the desired length. 

 We first recall the definition of c.m.~functions. 
\begin{Def}[Completely monotonic function \cite{widder1931necessary}]
	A function $f$ is said to be c.m.~on an interval $I$ if $f$ is
	continuous on $I$ and has derivatives of all orders on $\interior I$ (the interior of $I$) and for all $n \in \mathbb{N}$, 
	\begin{align}
		(-1)^n f^{(n)}(x) \geq 0, \quad x\in \interior I.
	\end{align}
\end{Def}

By \cite{widder1931necessary}, the sequence $(f(k))_{k=0}^\infty$ is c.m.~if $f$ is c.m.~on $[0,\infty)$. Thus, we can use c.m.~functions to generate moment sequences.

Based on the six types of moment sequences discussed in \Cref{sec:rateofdecay}, we define six classes of c.m.~functions that decay at different rates.  
\begin{Def}[Six classes of c.m.~functions]\label{def:classofcm}
	\leavevmode 
	Based on the definition of rate of decay in \Cref{def:rateofdecay}, we define the following c.m.~functions. 
	\begin{itemize}
		\item Sub-power-law decay $\Omega(n^{-s})$: $F_1(n; a, s) \triangleq a^s (\log(n + 1) + a)^{-s}$.
		\item Power-law decay $\Theta(n^{-s})$: $F_2(n;a,s) \triangleq a^s (n+a)^{-s} $.
		\item Soft power-law decay: $F_3 (n; a, b, s) \triangleq a^s b (n+a)^{-s}  (\log(n+1) + b)^{-1}$. 
		\item Intermediate decay: $F_4(n; a, s) \triangleq \exp( s \sqrt{a} - s \sqrt{n+a})$.
		\item Soft exponential decay: $F_5(n; a, s) \triangleq \exp(-sn)a (n+a)^{-1} $.
		\item Exponential decay $\Theta(\exp(-s n ))$: $F_6(n; s) \triangleq \exp( -sn)$.
	\end{itemize}
	$s \in \mathbb{R}^+$ is referred to as the rate, and $a\in \mathbb{R}^+$ and $b \in \mathbb{R}^+$ (if present) are referred to as the other parameters. 
\end{Def}

For moment sequences $(m_{i,n})_{n = 0}^{\infty}$, $i = 1,..., N$, the convex combination $m_n = \sum_{ i = 1}^N c_i m_{i,n}$, $\sum_{ i = 1}^N c_i = 1$, $c_i >0$, is still c.m., and thus $(m_n)_{n = 0}^{\infty}$ is also a moment sequence. 

Now, we introduce our algorithm to generate $(m_n)_{n = 0}^k$ with a desired rate of decay $s$ and flexibility in the other parameters, i.e., $a$ and/or $b$. {Since the c.m.~functions are analytic, the c.m.~function of integer moments can be easily extended to complex moments by analytic continuation. To obtain complex moments $m_z$, we just need to replace $n$ in steps $2$ and $3$ of \Cref{sch:cm} with $z \in \mathbb{C}$.}

\begin{algorithm}[H]
	\caption{Algorithm for moment sequences by c.m.~functions}\label{sch:cm}
	\begin{algorithmic}[1]
		\renewcommand{\algorithmicrequire}{\textbf{Input:} }
		\renewcommand{\algorithmicensure}{\textbf{Output:} }
		\REQUIRE A desired rate of decay (including the type of decay and the rate $s$), the number of mixed decays $N$, the other rates $s_i, 2\leq i \leq N$, the other parameters $a_i, b_i, \,1\leq i \leq N$, the weights $c_i, \,1\leq i \leq N$, and the number of moments $k$. 
		\ENSURE  $(m_n)_{n = 0}^k$   
		\STATE Based on the type of decay, set $F_0$ to be one of $F_1, F_2, F_3, F_4, F_5, F_6$ in \Cref{def:classofcm}. Suppose $F_0$ has three parameters $a_i, b_i, s$. If it is not the case, we will ignore $a_i$ and/or $b_i$ in the following steps. 
		\STATE Set $m_{i,n} = F_0(n; a_i, b_i, s_i), \, i \in [N], n \in [k]_0$. 
		\STATE Set $m_n = \sum_{ i = 1}^N c_i m_{i,n}, n \in [k]_0$. 
	\end{algorithmic}
\end{algorithm}

{To deterministically generate moments by \Cref{sch:cm}, we need to fix $s > 0, N  \in \mathbb{N}, s_i >s,  2\leq i \leq N$ for the exponential decay and $s_i = s$ for the others, $a_i, b_i, c_i > 0$ and $\sum_{ i = 1}^N c_i = 1$.  }

{To randomly generate moments by \Cref{sch:cm}, we first choose $s\sim \mathrm{Uniform}(0,10)$ except for the exponential and soft exponential cases, where $s\sim -\log\left(\mathrm{Uniform}(0,1)\right)$. We set $s_i = s$ except that $ s_1 = s, s_i \sim -\log\left( \mathrm{Uniform} (0,\exp(-s))\right),  2\leq i \leq N$ for the exponential decay. 
	Next, we randomly generate $N$ from the set $\{1,2,..., 10\}$. Then we choose $a_i, b_i  \sim \mathrm{Uniform} (0,10), \,1\leq i \leq N$. Finally we choose $c_i \sim \mathrm{Uniform}(0,10)$ and normalize them so that the sum is $1$.} 
In this way, we have at least $2N$ degrees of freedom.

Similar to \Cref{sec:cano}, the exact cdf is unknown. We use the GP approximations as a reference when making comparisons. 

\subsection{Generating moments by known distributions}\label{sec:bypdfs}
Beta distributions are continuous probability distributions supported on $[0,1]$. Here, we use the beta mixture model to randomly generate distributions. The pdf of $N$ beta mixtures is 
\begin{align}
	f_{\rm{beta}}(x,\mathbf{a, b, c}) \triangleq \sum_{ i = 1}^N c_i x^{a_i - 1}(1-x)^{b_i - 1}, 
\end{align}
where $\mathbf{a} \triangleq (a_i)_{i = 1}^N, \mathbf{b} \triangleq (b_i)_{i = 1}^N, \mathbf{c} \triangleq (c_i)_{i = 1}^N$, $\mathbf{a, b, c} \succ \mathbf{0}$ and $\sum_{ i = 1}^N c_i = 1$. The $n$-th moment is
\begin{align}
	m_n = \sum_{ i = 1}^N c_i\prod_{r=0}^{n-1} \frac{a_i+r}{a_i+b_i+r}, 
\end{align}
and the cdf is 
\begin{align}
	F_{\rm{beta}}(x,\mathbf{a, b, c}) = \sum_{ i = 1}^N c_i \int_{0}^{x}t^{a_i - 1}(1-t)^{b_i - 1}\, dt. 
\end{align}

By randomly choosing $\mathbf{a, b, c} $, we are able to randomly generate the distribution and thus randomly generate the imaginary and integer moment sequences by definition. We first randomly generate $N$ from the set $\{1,2,..., 10\}$.\footnote{In order to increase the likelihood of generating moment sequences that can span the entire space, we randomly generate $N$, instead of a fixed value. Specifically, for a small number of realizations, such as $100$, the span behavior is observed to vary with different values of $N$. Thus, randomizing $N$ serves to increase the likelihood of achieving the desired outcome. } Then we choose $a_i, b_i  \sim \mathrm{Uniform} (0,10), \,1\leq i \leq N$. Finally we choose $c_i \sim \mathrm{Uniform}(0,10)$ and normalize them so that the sum is $1$.

\subsection{Performance evaluation}\label{sec:perf}

Using the moments randomly generated in \Cref{sec:momentsbyanalytic}, we compare the performance of the different Hausdorff moment transforms.  
Since the actual cdfs are not available in this case, we use the cdfs obtained by the GP method as the reference.\footnote{As the GP method utilizes imaginary moments, it naturally processes  the complete information to recover the exact distribution. However, the HMTs that produce solutions to the THMP, such as the ME method, may not be close to the exact one as they only utilize a finite number of integer moments. The gap does not necessarily indicate that an HMT performs poorly—the method may simply assume a different continuation of the moments sequence. }

\begin{figure}
	\centering
	\subfloat[Sub power-law decay]{\includegraphics[width=0.31\figwidth]{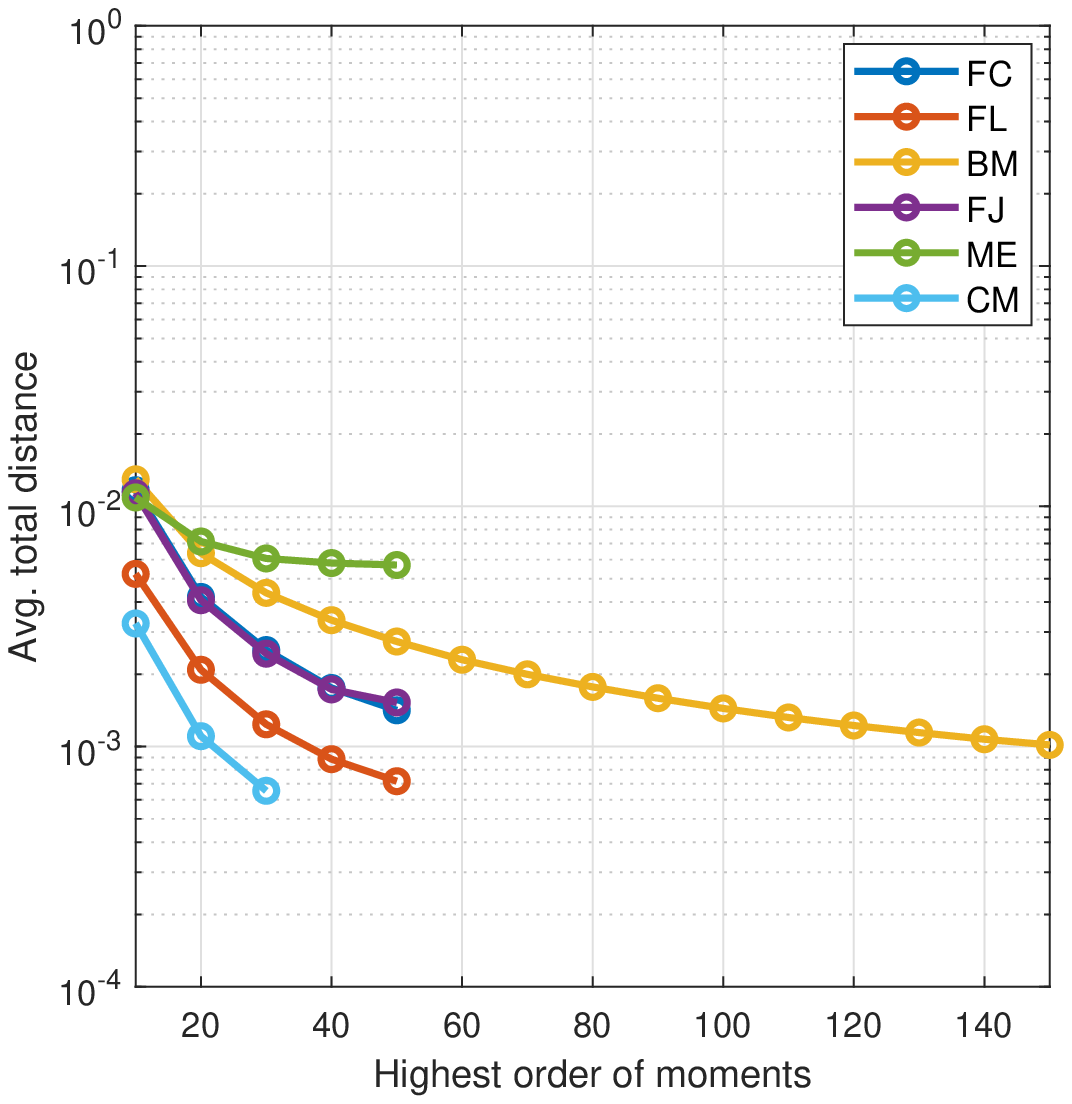}}\quad
	\subfloat[Power-law decay\label{fig:pow}]{\includegraphics[width=0.31\figwidth]{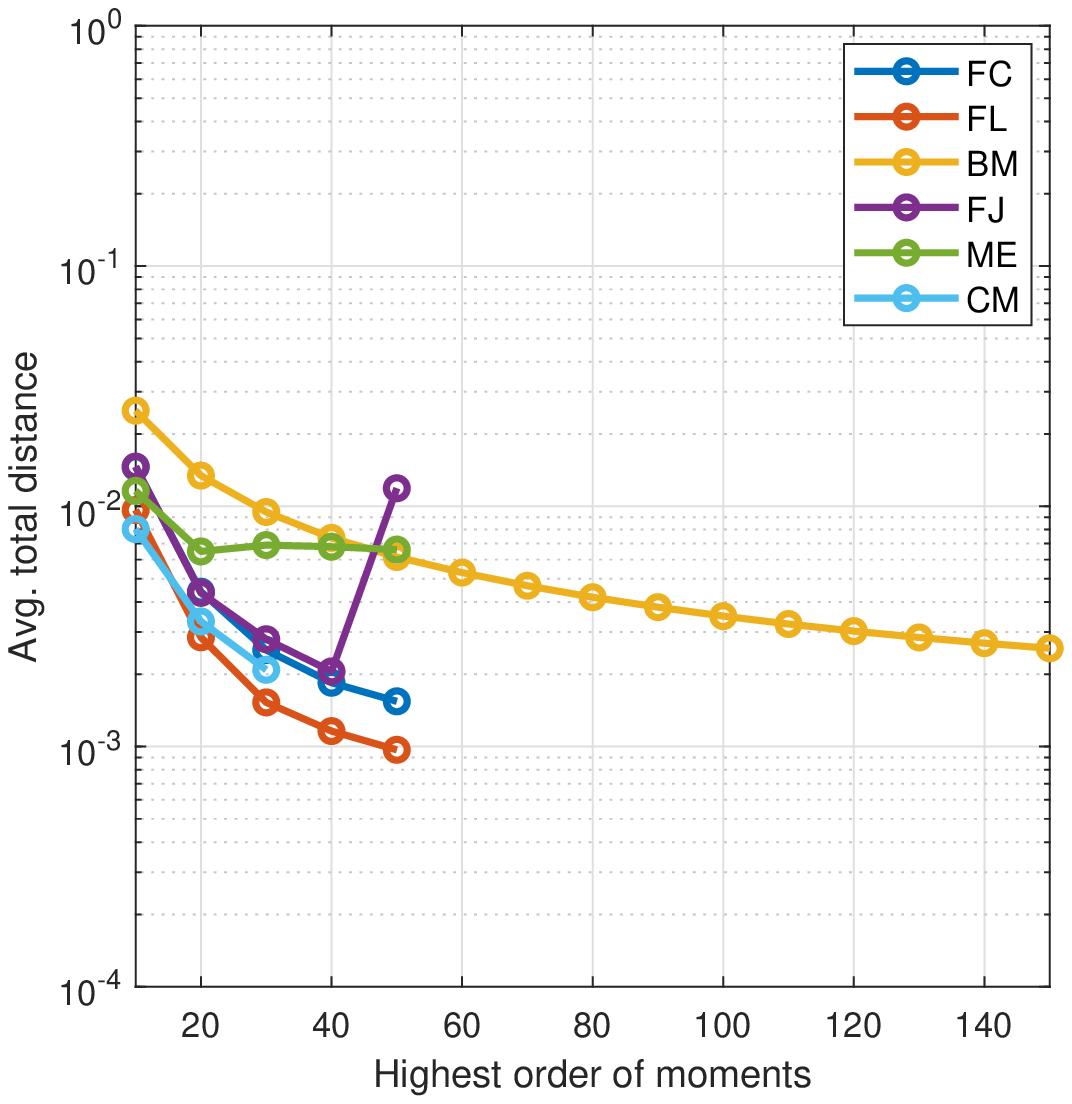}}\quad
	\subfloat[Soft power-law decay]{\includegraphics[width=0.31\figwidth]{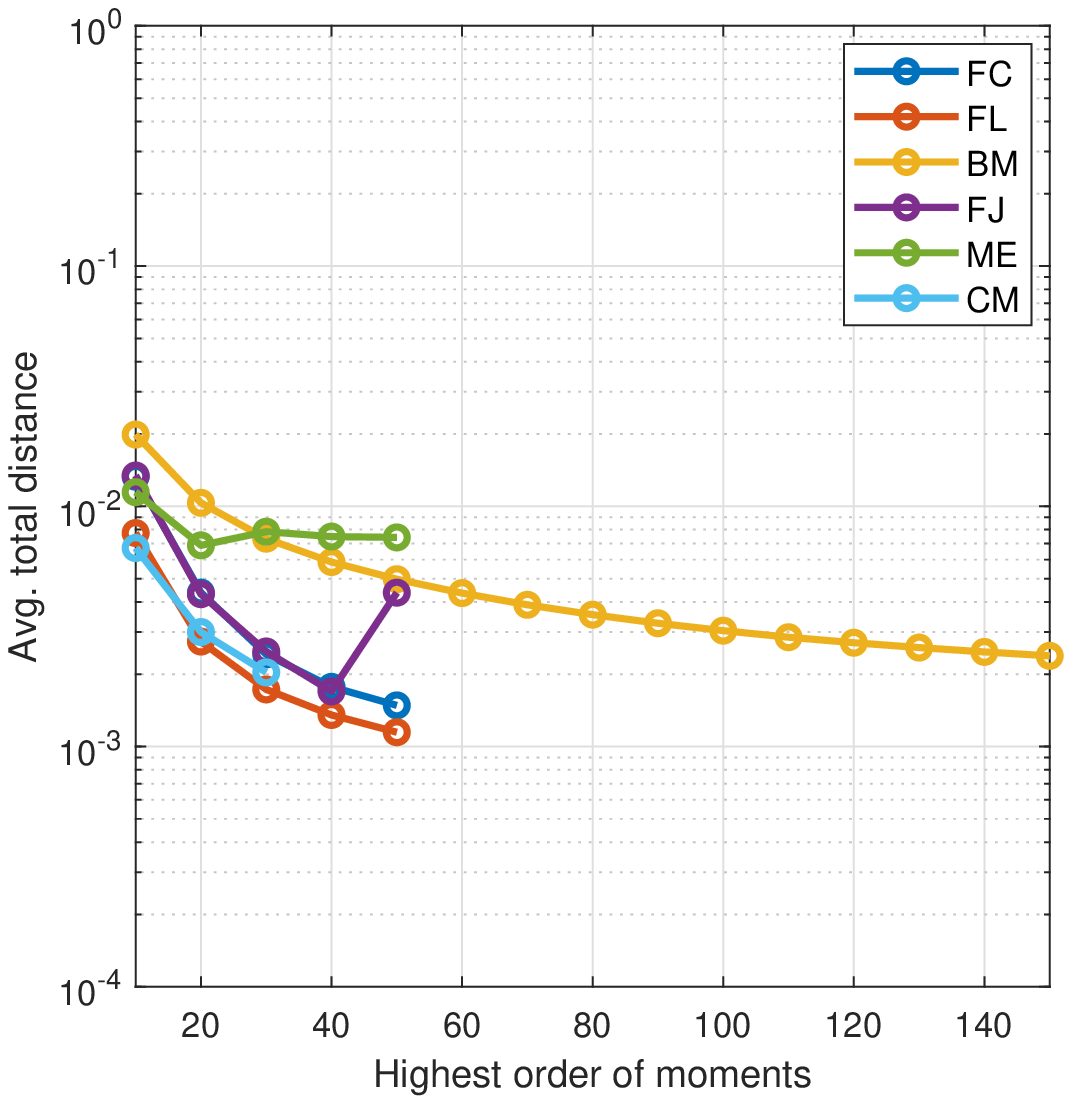}}\\
	\subfloat[Intermediate decay]{\includegraphics[width=0.31\figwidth]{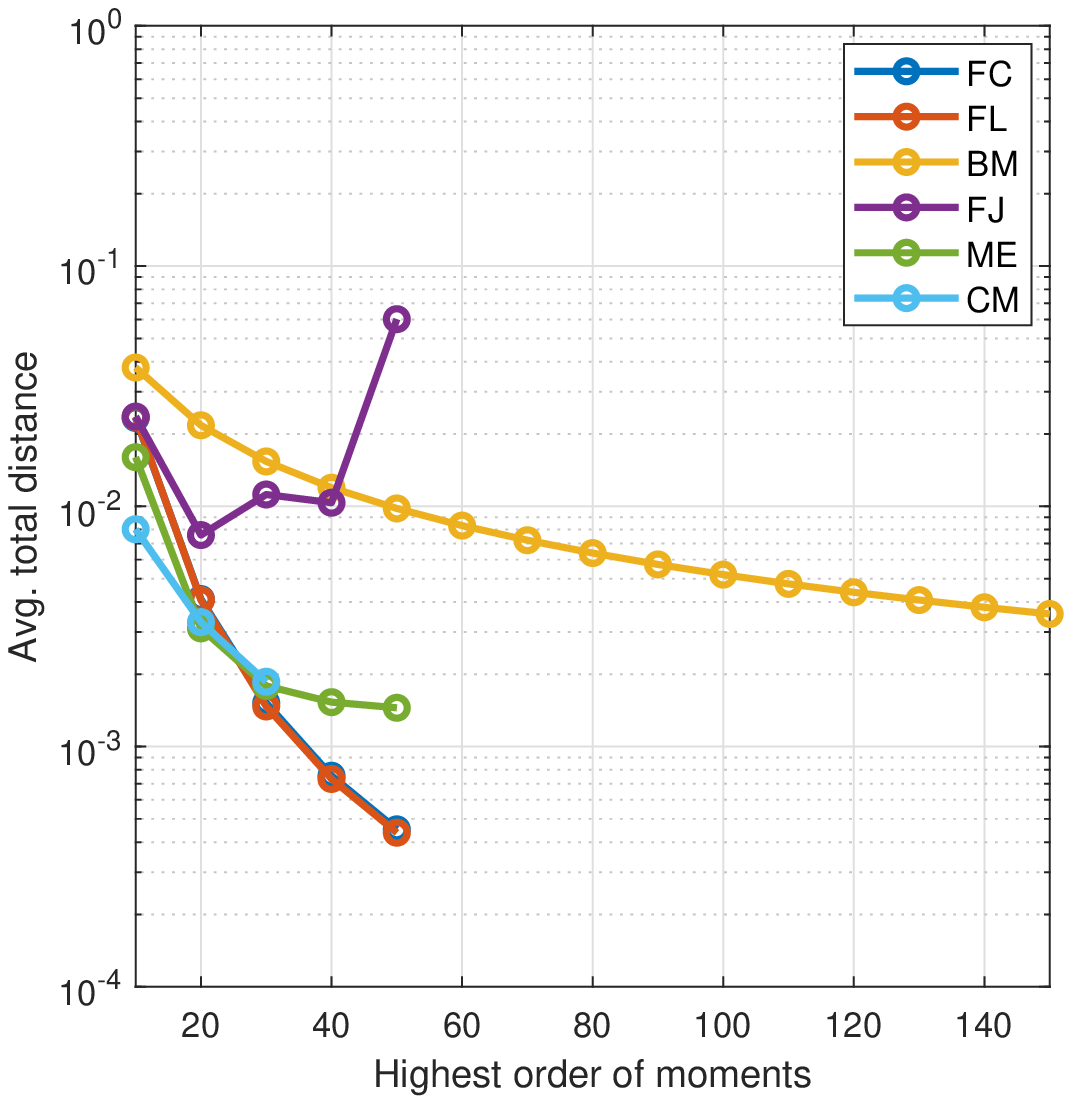}}\quad
	\subfloat[Exponential decay \label{fig:dvvsnExp}]{\includegraphics[width=0.31\figwidth]{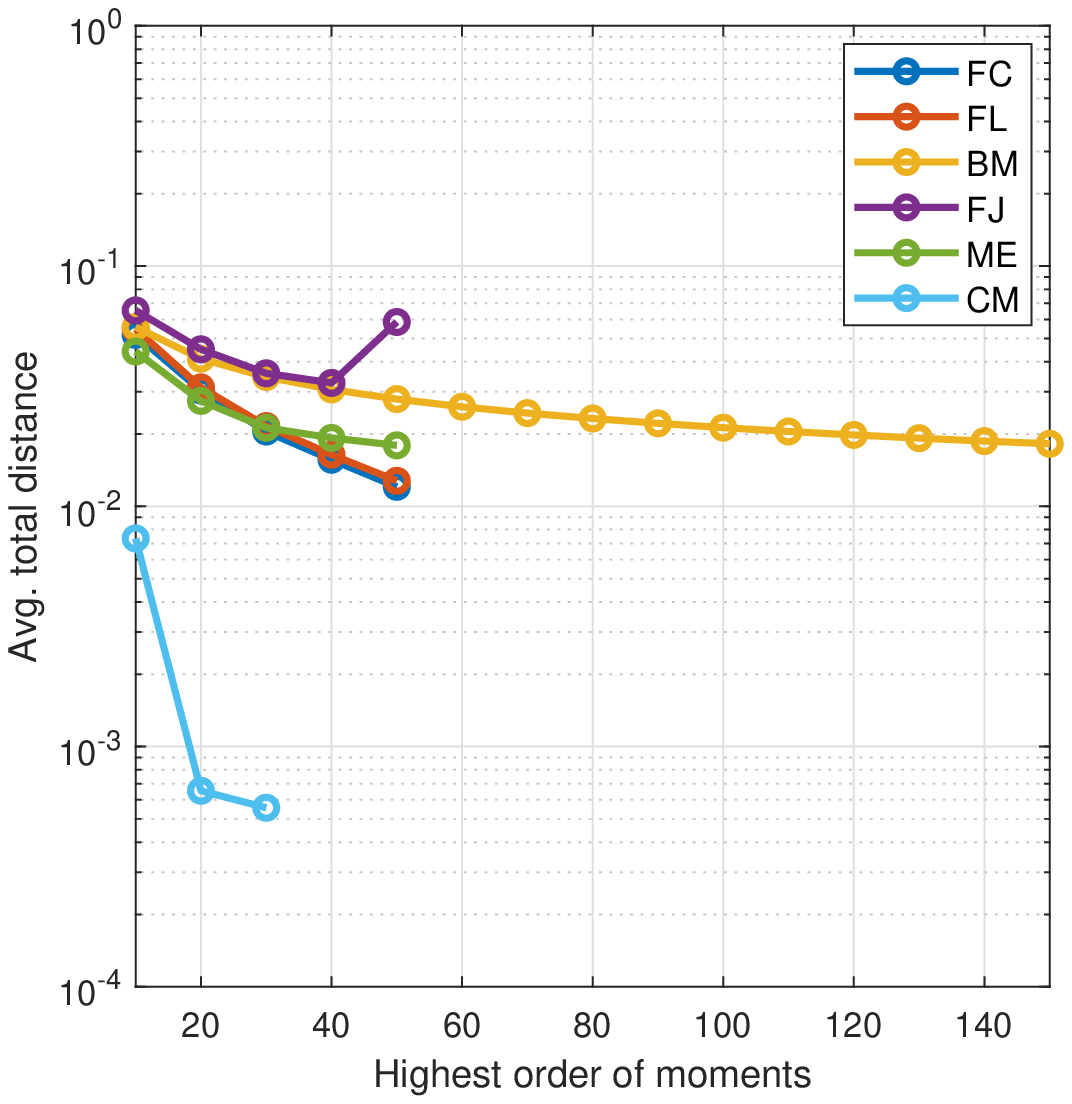}}\quad
	\subfloat[Soft exponential decay]{\includegraphics[width=0.31\figwidth]{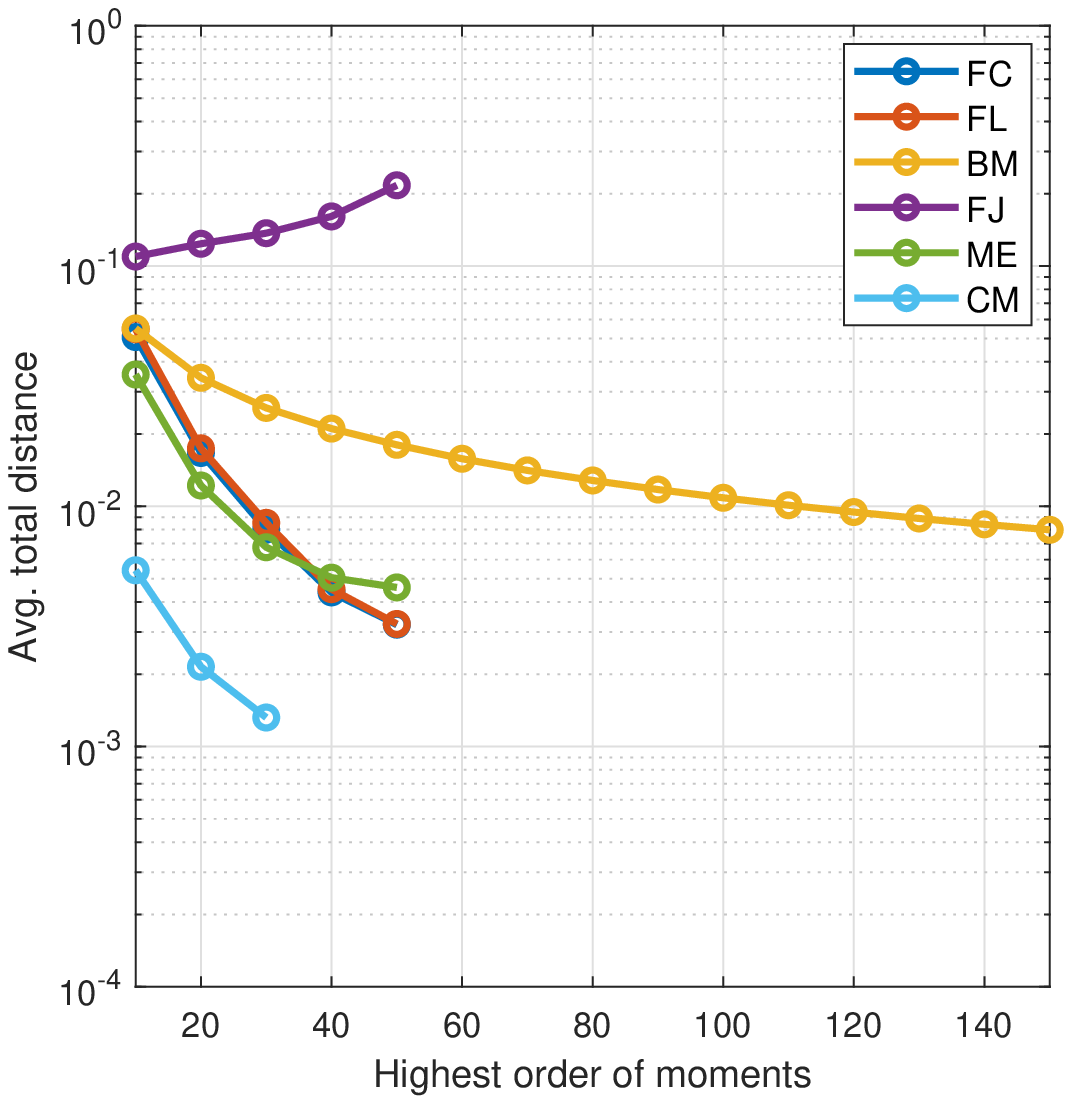}}
	\caption{The average of the total distance versus the number of moments for different methods and different types of decays. Averaging is performed over $100$ randomly generated moment sequences. }\label{fig:dvvsn}
\end{figure}

\begin{figure}
	\centering
	\subfloat[Original implementation.\label{fig:comptimeoriginal}]{\includegraphics[width=0.45\figwidth]{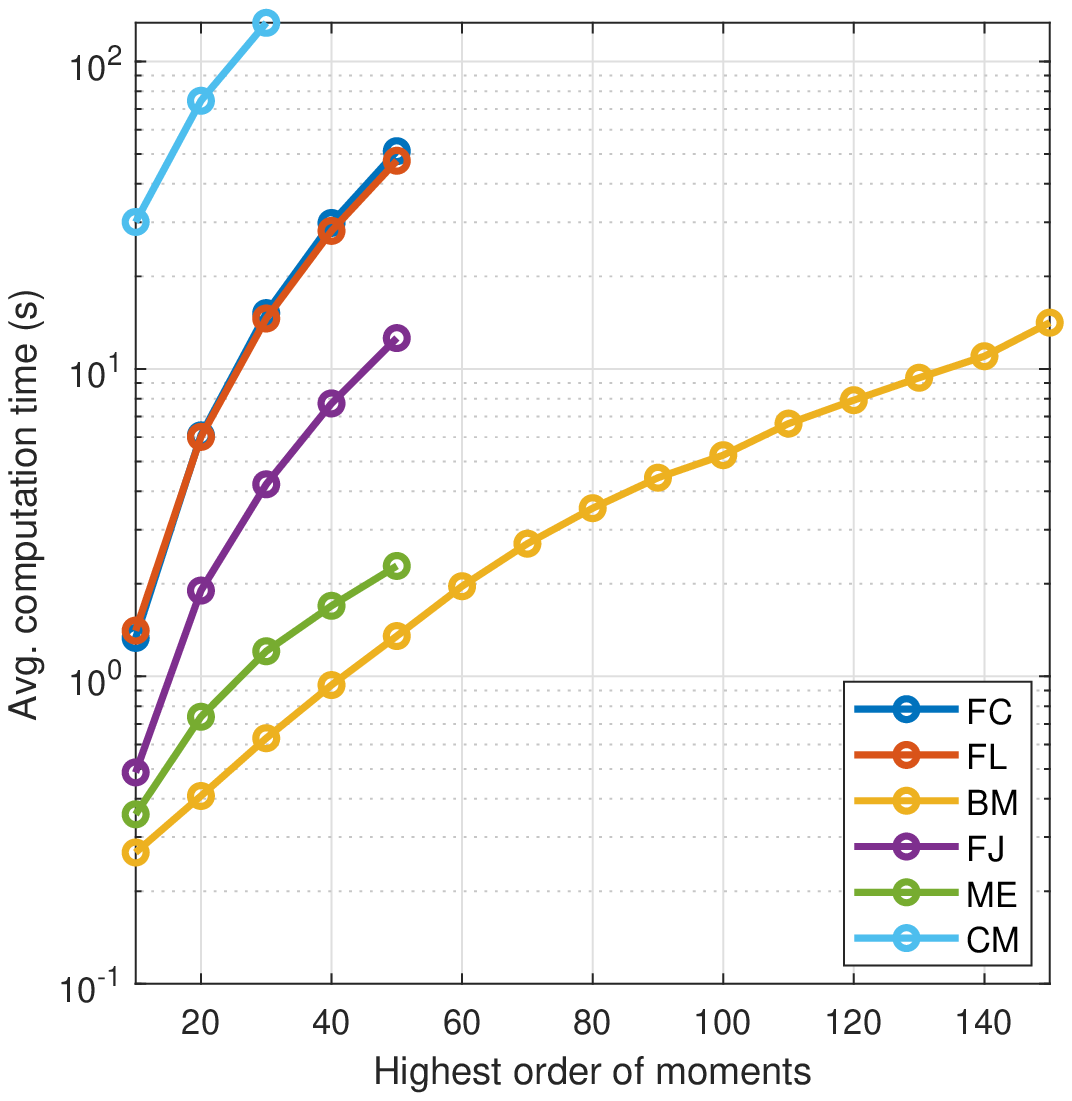}}\quad
	\subfloat[Linear transform. \label{fig:comptimelinear}]{\includegraphics[width=0.45\figwidth]{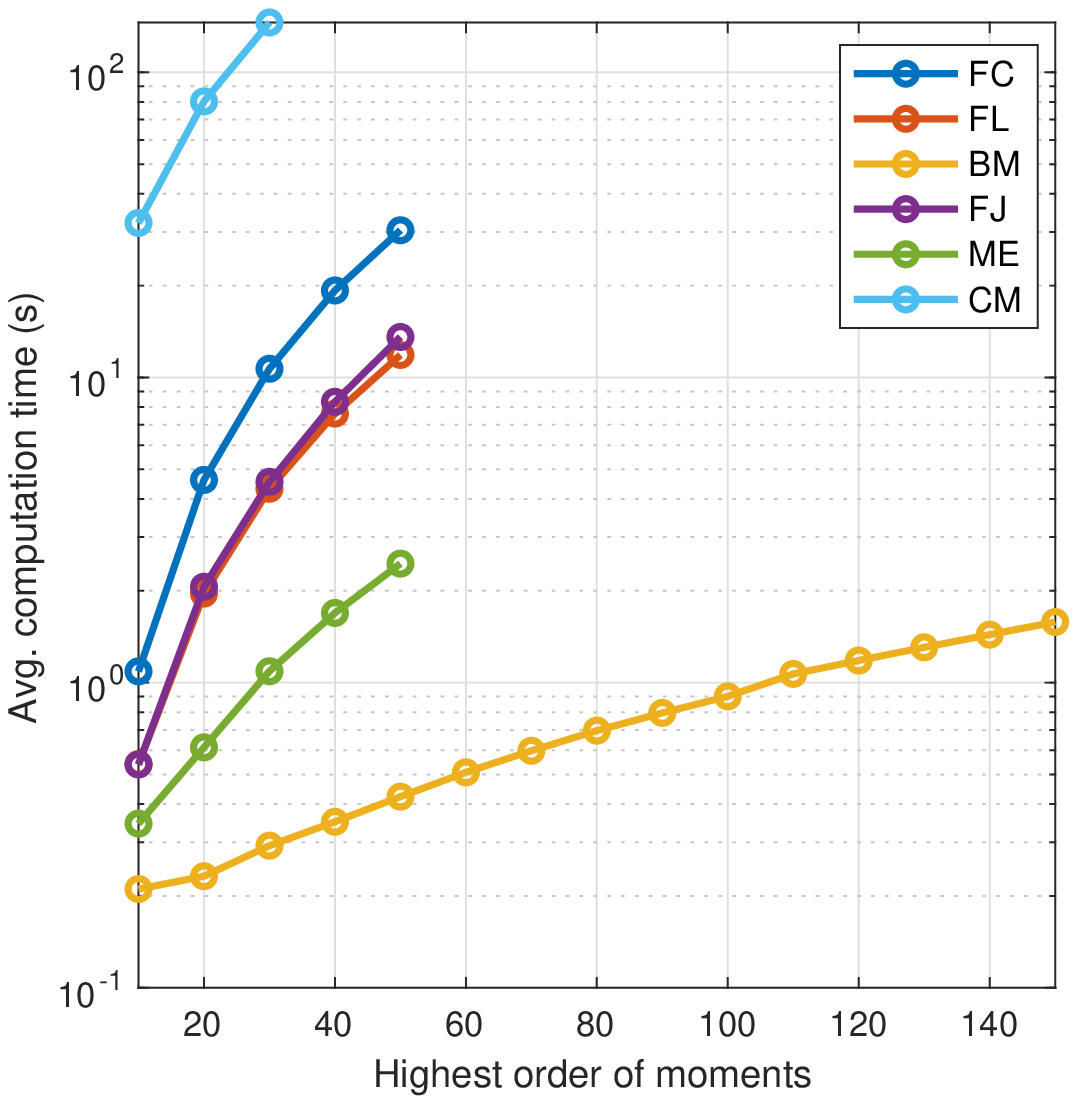}
	} 
	\caption{The average computation time versus the number of moments. Averaging is performed over $100$ randomly generated moment sequences of power-law decay.}\label{fig:powerlaw} 
\end{figure}

For all methods except the GP, BM, FL, and FC methods, we consider  $F_{\rm{method},n}|_{\mathcal{ U}_n}$. For the GP method, we consider $F_{\rm{GP}}|_{\mathcal{ U}_{150}}$, for the BM method, we consider $F_{\rm{BM},n}|_{\mathcal{ U}_{n+1}}$, and for the FL and FC methods, we consider $F_{\rm{method},n-1}|_{\mathcal{ U}_{n}}$. We apply the tweaking mapping in \Cref{def:correction} followed by monotone cubic interpolation. 

In \Cref{fig:dvvsn}, we plot the average of the total distance versus the number of moments ($m_0$ is not counted) for different methods and different types of decays. For practical purposes regarding accuracy requirements in moments and basic operations and the difficulty in finding roots of polynomials, the number of moments of the FJ, FL, FC, and ME methods are limited to $50$, the number of moments of the CM method is limited to $30$, and the number of moments of the BM method is limited to $150$. {We observe that the average of the total distance of all the methods except the FJ and the ME methods decreases as the number of moments increases. In most cases, the average of the total distance of the FJ method does not decrease since it is not guaranteed to converge. As for the ME method, although it converges in entropy, it suffers from the problem of finding the minimum in \eqref{eq:me_cvx}. Even if some of the curves in \Cref{fig:dvvsn} show a decreasing trend, we cannot rely on this method as long as we cannot make sure the results we get from solving \eqref{eq:me_cvx} are the minima, i.e., the size of the gradient is less than
	the value of the optimality tolerance.
	In terms of their best possible performance, we also observe that the CM method is great for all types of decays, especially the exponential and soft exponential decays. For the exponential decay, the moments are from discrete distributions that have no more than $10$ jumps. {Thus, the raw CM method of order $20$ is able to fully characterize the distributions by identifying jump locations and jump heights.} Our comparisons are made based on the decay type. However, in practice, it is impossible to determine the decay type given the values of the moments. Thus, for practical purposes, it is preferable to utilize methods that work well overall. Among all the methods, the FL and FC methods have overall good performance and the FL method even outperforms the CM method in the power-law decay. }

\Cref{fig:powerlaw} shows computation time versus the number of moments for the power-law decay and it is similar for the other types. In \Cref{fig:comptimeoriginal}, the computation time of the FL, FC, and BM method is based on the original implementation, whereas in \Cref{fig:comptimelinear}, it is based on the linear transform with the pre-calculated transform matrices.  In \Cref{fig:comptimelinear}, for all methods except the CM method, the computation time increases quadratically with the number of moments. This is consistent with the order of matrix-vector multiplication for the FL, FC and BM method. Furthermore, with the linear transform, the computation time for the FL, FC and BM method in \Cref{fig:comptimeoriginal} is significantly reduced compared to that in \Cref{fig:comptimeoriginal}. Even for small values of $n$, such as $10$, $20$, and $30$, the CM method's computation time is at least ten times longer than that of the other methods for the same number of moments. As $n$ increases, the computation time of the CM method is expected to grow even more, making it unsuitable for time-sensitive applications. 

{As shown in \Cref{exam:beta}, the moment sequences generated by \Cref{sec:bypdfs} exhibit power-law decay. Its results are similar to those in \Cref{fig:pow}.}

\subsection{Accuracy requirements for HMTs}\label{sec:accuracyHMT}
The accuracy requirements apply to both of the coefficients of the moments and the HMTs, such as the transform matrices for the BM and FL method. Here we just mention some typical requirements.  

Accuracy requirements for transform matrices: For the BM method of order $n$, as $\max \|\mathbf{A}\| \sim \frac{\sqrt{27}}{2\pi} \frac{3^n}{n}$, $n \to \infty$, about $n/2 - \log_{10} n$ decimal digits are needed \cite{haenggi2018efficient}. For the FL method of order $n$, as $\max \|\mathbf{\hat{A}}\| \leq \binom{2n}{n}\binom{n}{n/3} \sim \frac{3^{n+1}2^{4n/3}}{4\pi n}$, $n \to \infty$, about $0.9n $ decimal digits are sufficient. 

Accuracy requirements for moments in \Cref{def:classofcm}: For the power-law decay, about $s \log_{10} \frac{n+a}{a}$ decimal digits are needed. For exponential decay, about $sn/2$ decimal digits are needed. 

The calculation can be carried out in the logarithm domain which may reduce the accuracy requirement.  Moreover, the matrix multiplication of the FL method is pre-calculated, i.e., $\hat{\mathbf{A}} \mathbf{1}$ and $\hat{\mathbf{A}} \hat{\mathbf{m}}$ are done separately.

\section{Conclusion}
The truncated Hausdorff moment problem is both fundamental and practically relevant. While several theoretical methods for its solution have been proposed, a systematic study and comparison of their performance and details in their numerical implementation has been missing. This work provides these important missing pieces by discussing advantages and drawbacks of different methods, comparing convergence properties, and proposing a tweaking method to ensure the resulting functions satisfy the properties of probability distributions. In addition, two new methods are introduced that compare favorably with known ones. For a fair comparison, three methods for the random generation of moment sequences are discussed or proposed, so that the distances between approximated and actual distributions can be averaged over many realizations of moment sequences. It turns out that the decrease of the distances as a function of the number of moments depends strongly on the type of decay of the moment sequences. For this reason, we introduced a classification of the decay order in six classes. Our work also reveals the trade-off between the performance and computational cost. If high accuracy is desired (at higher computational effort), the FC and FL methods are good choices, with the FC method having the advantage of uniform convergence; if computational resources are scarce, the BM method is preferred.

\bibliography{ref}

\end{document}